\documentclass[11pt,a4paper]{article}

\usepackage{slashed}
\usepackage{xkeyval}
\usepackage{url,amsmath,amssymb,latexsym,pstricks,mathrsfs,comment,amsthm,graphicx,tikz,tikz-cd,enumerate,accents,pgffor,cite,wrapfig,multicol,float}
\usepackage{bibspacing}
\newcommand{\nc}{\newcommand}
\nc{\rnc}{\renewcommand}

\nc{\OEIS}{}

\usepackage{geometry}  \rnc{\ss}{\smallskip} \nc{\ms}{\medskip}  \nc{\nss}{\vspace{-3mm}} 

\nc{\uvertcols}[2]{\foreach \x in {#1}{ \uvertcol{\x}{#2}}}
\nc{\lvertcols}[2]{\foreach \x in {#1}{ \lvertcol{\x}{#2}}}
\nc{\uverts}[1]{\foreach \x in {#1}{ \uvert{\x}}}
\nc{\lverts}[1]{\foreach \x in {#1}{ \lvert{\x}}}
\nc{\Bnt}{\B_n^\tau}
\nc{\DClass}[2]{D_{#1;#2}}
\nc{\Ideal}[2]{I_{#1;#2}}
\nc{\Minimal}[2]{M_{#1;#2}}
\nc{\leqt}{\leq^{\tau}}
\nc{\DC}[3]{
\draw[thick,fill=lightgray](#1-2.5,#2-1)--(#1+2.5,#2-1)--(#1+2.5,#2+1)--(#1-2.5,#2+1)--(#1-2.5,#2-1);
\draw(#1,#2)node{$#3$};
}
\nc{\DColumn}[3]{
\draw(#1,#2)--(#1,#2+6);
\draw 
(#1,#2+6) node [rounded corners,rectangle,draw,fill=blue!20] {$\DClass{7}{#3}$}
(#1,#2+4) node [rounded corners,rectangle,draw,fill=blue!20] {$\DClass{5}{#3}$}
(#1,#2+2) node [rounded corners,rectangle,draw,fill=blue!20] {$\DClass{3}{#3}$}
(#1,#2+0) node [rounded corners,rectangle,draw,fill=blue!20] {$\DClass{1}{#3}$};
}
\nc{\DRow}[1]{\draw(0,#1)--(12+1.5,#1-4.5);\draw[dashed](12+1.5,#1-4.5)--(12+3,#1-5);}
\nc{\dla}{\la\!\la}
\nc{\dra}{\ra\!\ra}
\nc{\BnSn}{\B_n\sm\S_n}
\nc{\BntSn}{\Bnt\sm\S_n}
\nc{\Et}{E^\tau}
\nc{\bbEt}{\bbE^\tau}

\nc{\DCo}[3]{
\draw(#1,#2)--(#1,#2+10);
\draw 
(#1,#2+10) node [rounded corners,rectangle,draw,fill=blue!20] {$\DClass{11}{#3}$}
(#1,#2+8) node [rounded corners,rectangle,draw,fill=blue!20] {$\DClass{9}{#3}$}
(#1,#2+6) node [rounded corners,rectangle,draw,fill=blue!20] {$\DClass{7}{#3}$}
(#1,#2+4) node [rounded corners,rectangle,draw,fill=blue!20] {$\DClass{5}{#3}$}
(#1,#2+2) node [rounded corners,rectangle,draw,fill=blue!20] {$\DClass{3}{#3}$}
(#1,#2+0) node [rounded corners,rectangle,draw,fill=blue!20] {$\DClass{1}{#3}$};
}
\nc{\DRo}[1]{\draw(0,#1)--(18+1.5,#1-6.5);\draw[dashed](18+1.5,#1-6.5)--(18+3,#1-7);}

\makeatletter

\DeclareMathSymbol{\widehatsym}{\mathord}{largesymbols}{"62}
\newcommand\lowerwidehatsym{%
  \text{\smash{\raisebox{-1.3ex}{%
    $\widehatsym$}}}}
\newcommand\fixwidehat[1]{%
  \mathchoice
    {\accentset{\displaystyle\lowerwidehatsym}{#1}}
    {\accentset{\textstyle\lowerwidehatsym}{#1}}
    {\accentset{\scriptstyle\lowerwidehatsym}{#1}}
    {\accentset{\scriptscriptstyle\lowerwidehatsym}{#1}}
}
\rnc{\widehat}{\fixwidehat}

\begin{document}


\nc{\ubluebox}[2]{\bluebox{#1}{1.7}{#2}2\udotted{#1}{#2}}
\nc{\lbluebox}[2]{\bluebox{#1}0{#2}{.3}\ldotted{#1}{#2}}
\nc{\ublueboxes}[1]{{
\foreach \x/\y in {#1}
{ \ubluebox{\x}{\y}}}
}
\nc{\lblueboxes}[1]{{
\foreach \x/\y in {#1}
{ \lbluebox{\x}{\y}}}
}

\nc{\bluebox}[4]{
\draw[color=blue!20, fill=blue!20] (#1,#2)--(#3,#2)--(#3,#4)--(#1,#4)--(#1,#2);
}
\nc{\redbox}[4]{
\draw[color=red!20, fill=red!20] (#1,#2)--(#3,#2)--(#3,#4)--(#1,#4)--(#1,#2);
}

\nc{\rectangle}[4]{
\draw (#1,#2)--(#3,#2)--(#3,#4)--(#1,#4)--(#1,#2);
}

\nc{\bluetrap}[8]{
\draw[color=blue!20, fill=blue!20] (#1,#2)--(#3,#4)--(#5,#6)--(#7,#8)--(#1,#2);
}
\nc{\redtrap}[8]{
\draw[color=red!20, fill=red!20] (#1,#2)--(#3,#4)--(#5,#6)--(#7,#8)--(#1,#2);
}

\usetikzlibrary{decorations.markings}
\usetikzlibrary{arrows,matrix}
\usepgflibrary{arrows}
\tikzset{->-/.style={decoration={
  markings,
  mark=at position #1 with {\arrow{>}}},postaction={decorate}}}
\tikzset{-<-/.style={decoration={
  markings,
  mark=at position #1 with {\arrow{<}}},postaction={decorate}}}
\nc{\Unode}[1]{\draw(#1,-2)node{$U$};}
\nc{\Dnode}[1]{\draw(#1,-2)node{$D$};}
\nc{\Fnode}[1]{\draw(#1,-2)node{$F$};}
\nc{\Cnode}[1]{\draw(#1-.1,-2)node{$\phantom{+0},$};}
\nc{\Unodes}[1]{\foreach \x in {#1}{ \Unode{\x} }}
\nc{\Dnodes}[1]{\foreach \x in {#1}{ \Dnode{\x} }}
\nc{\Fnodes}[1]{\foreach \x in {#1}{ \Fnode{\x} }}
\nc{\Cnodes}[1]{\foreach \x in {#1}{ \Cnode{\x} }}
\nc{\Uedge}[2]{\draw[->-=0.6,line width=.3mm](#1,#2-9)--(#1+1,#2+1-9); \vertsm{#1}{#2-9} \vertsm{#1+1}{#2+1-9}}
\nc{\Dedge}[2]{\draw[->-=0.6,line width=.3mm](#1,#2-9)--(#1+1,#2-1-9); \vertsm{#1}{#2-9} \vertsm{#1+1}{#2-1-9}}
\nc{\Fedge}[2]{\draw[->-=0.6,line width=.3mm](#1,#2-9)--(#1+1,#2-9); \vertsm{#1}{#2-9} \vertsm{#1+1}{#2-9}}
\nc{\Uedges}[1]{\foreach \x/\y in {#1}{\Uedge{\x}{\y}}}
\nc{\Dedges}[1]{\foreach \x/\y in {#1}{\Dedge{\x}{\y}}}
\nc{\Fedges}[1]{\foreach \x/\y in {#1}{\Fedge{\x}{\y}}}
\nc{\xvertlabel}[1]{\draw(#1,-10+.6)node{{\tiny $#1$}};}
\nc{\yvertlabel}[1]{\draw(0-.4,-9+#1)node{{\tiny $#1$}};}
\nc{\xvertlabels}[1]{\foreach \x in {#1}{ \xvertlabel{\x} }}
\nc{\yvertlabels}[1]{\foreach \x in {#1}{ \yvertlabel{\x} }}

\nc{\bbE}{\mathbb E}
\nc{\floorn}{\lfloor\tfrac n2\rfloor}
\rnc{\sp}{\supseteq}
\rnc{\arraystretch}{1.2}

\nc{\bn}{{[n]}} \nc{\bt}{{[t]}} \nc{\ba}{{[a]}} \nc{\bl}{{[l]}} \nc{\bm}{{[m]}} \nc{\bk}{{[k]}} \nc{\br}{{[r]}} \nc{\bs}{{[s]}} \nc{\bnf}{{[n-1]}}

\nc{\M}{\mathcal M}
\nc{\G}{\mathcal G}
\nc{\F}{\mathbb F}
\nc{\MnJ}{\mathcal M_n^J}
\nc{\EnJ}{\mathcal E_n^J}
\nc{\Mat}{\operatorname{Mat}}
\nc{\RegMnJ}{\Reg(\MnJ)}
\nc{\row}{\mathfrak r}
\nc{\col}{\mathfrak c}
\nc{\Row}{\operatorname{Row}}
\nc{\Col}{\operatorname{Col}}
\nc{\Span}{\operatorname{span}}
\nc{\mat}[4]{\left[\begin{matrix}#1&#2\\#3&#4\end{matrix}\right]}
\nc{\tmat}[4]{\left[\begin{smallmatrix}#1&#2\\#3&#4\end{smallmatrix}\right]}
\nc{\ttmat}[4]{{\tiny \left[\begin{smallmatrix}#1&#2\\#3&#4\end{smallmatrix}\right]}}
\nc{\tmatt}[9]{\left[\begin{smallmatrix}#1&#2&#3\\#4&#5&#6\\#7&#8&#9\end{smallmatrix}\right]}
\nc{\ttmatt}[9]{{\tiny \left[\begin{smallmatrix}#1&#2&#3\\#4&#5&#6\\#7&#8&#9\end{smallmatrix}\right]}}
\nc{\MnGn}{\M_n\sm\G_n}
\nc{\MrGr}{\M_r\sm\G_r}
\nc{\qbin}[2]{\left[\begin{matrix}#1\\#2\end{matrix}\right]_q}
\nc{\tqbin}[2]{\left[\begin{smallmatrix}#1\\#2\end{smallmatrix}\right]_q}
\nc{\qbinx}[3]{\left[\begin{matrix}#1\\#2\end{matrix}\right]_{#3}}
\nc{\tqbinx}[3]{\left[\begin{smallmatrix}#1\\#2\end{smallmatrix}\right]_{#3}}
\nc{\MNJ}{\M_nJ}
\nc{\JMN}{J\M_n}
\nc{\RegMNJ}{\Reg(\MNJ)}
\nc{\RegJMN}{\Reg(\JMN)}
\nc{\RegMMNJ}{\Reg(\MMNJ)}
\nc{\RegJMMN}{\Reg(\JMMN)}
\nc{\Wb}{\overline{W}}
\nc{\Xb}{\overline{X}}
\nc{\Yb}{\overline{Y}}
\nc{\Zb}{\overline{Z}}
\nc{\Sib}{\overline{\Si}}
\nc{\Om}{\Omega}
\nc{\Omb}{\overline{\Om}}
\nc{\Gab}{\overline{\Ga}}
\nc{\qfact}[1]{[#1]_q!}
\nc{\smat}[2]{\left[\begin{matrix}#1&#2\end{matrix}\right]}
\nc{\tsmat}[2]{\left[\begin{smallmatrix}#1&#2\end{smallmatrix}\right]}
\nc{\hmat}[2]{\left[\begin{matrix}#1\\#2\end{matrix}\right]}
\nc{\thmat}[2]{\left[\begin{smallmatrix}#1\\#2\end{smallmatrix}\right]}
\nc{\LVW}{\mathcal L(V,W)}
\nc{\KVW}{\mathcal K(V,W)}
\nc{\LV}{\mathcal L(V)}
\nc{\RegLVW}{\Reg(\LVW)}
\nc{\sM}{\mathscr M}
\nc{\sN}{\mathscr N}
\rnc{\iff}{\ \Leftrightarrow\ }
\nc{\Hom}{\operatorname{Hom}}
\nc{\End}{\operatorname{End}}
\nc{\Aut}{\operatorname{Aut}}
\nc{\Lin}{\mathcal L}
\nc{\Hommn}{\Hom(V_m,V_n)}
\nc{\Homnm}{\Hom(V_n,V_m)}
\nc{\Homnl}{\Hom(V_n,V_l)}
\nc{\Homkm}{\Hom(V_k,V_m)}
\nc{\Endm}{\End(V_m)}
\nc{\Endn}{\End(V_n)}
\nc{\Endr}{\End(V_r)}
\nc{\Autm}{\Aut(V_m)}
\nc{\Autn}{\Aut(V_n)}
\nc{\MmnJ}{\M_{mn}^J}
\nc{\MmnA}{\M_{mn}^A}
\nc{\MmnB}{\M_{mn}^B}
\nc{\Mmn}{\M_{mn}}
\nc{\Mkl}{\M_{kl}}
\nc{\Mnm}{\M_{nm}}
\nc{\EmnJ}{\mathcal E_{mn}^J}
\nc{\MmGm}{\M_m\sm\G_m}
\nc{\RegMmnJ}{\Reg(\MmnJ)}
\rnc{\implies}{\ \Rightarrow\ }
\nc{\DMmn}[1]{D_{#1}(\Mmn)}
\nc{\DMmnJ}[1]{D_{#1}(\MmnJ)}
\nc{\MMNJ}{\Mmn J}
\nc{\JMMN}{J\Mmn}
\nc{\JMMNJ}{J\Mmn J}
\nc{\Inr}{\mathcal I(V_n,W_r)}
\nc{\Lnr}{\mathcal L(V_n,W_r)}
\nc{\Knr}{\mathcal K(V_n,W_r)}
\nc{\Imr}{\mathcal I(V_m,W_r)}
\nc{\Kmr}{\mathcal K(V_m,W_r)}
\nc{\Lmr}{\mathcal L(V_m,W_r)}
\nc{\Kmmr}{\mathcal K(V_m,W_{m-r})}
\nc{\tr}{{\operatorname{T}}}
\nc{\MMN}{\MmnA(\F_1)}
\nc{\MKL}{\Mkl^B(\F_2)}
\nc{\RegMMN}{\Reg(\MmnA(\F_1))}
\nc{\RegMKL}{\Reg(\Mkl^B(\F_2))}
\nc{\gRhA}{\widehat{\mathscr R}^A}
\nc{\gRhB}{\widehat{\mathscr R}^B}
\nc{\gLhA}{\widehat{\mathscr L}^A}
\nc{\gLhB}{\widehat{\mathscr L}^B}
\nc{\timplies}{\Rightarrow}
\nc{\tiff}{\Leftrightarrow}
\nc{\Sija}{S_{ij}^a}
\nc{\dmat}[8]{\draw(#1*1.5,#2)node{$\left[\begin{smallmatrix}#3&#4&#5\\#6&#7&#8\end{smallmatrix}\right]$};}
\nc{\bdmat}[8]{\draw(#1*1.5,#2)node{${\mathbf{\left[\begin{smallmatrix}#3&#4&#5\\#6&#7&#8\end{smallmatrix}\right]}}$};}
\nc{\rdmat}[8]{\draw(#1*1.5,#2)node{\rotatebox{90}{$\left[\begin{smallmatrix}#3&#4&#5\\#6&#7&#8\end{smallmatrix}\right]$}};}
\nc{\rldmat}[8]{\draw(#1*1.5-0.375,#2)node{\rotatebox{90}{$\left[\begin{smallmatrix}#3&#4&#5\\#6&#7&#8\end{smallmatrix}\right]$}};}
\nc{\rrdmat}[8]{\draw(#1*1.5+.375,#2)node{\rotatebox{90}{$\left[\begin{smallmatrix}#3&#4&#5\\#6&#7&#8\end{smallmatrix}\right]$}};}
\nc{\rfldmat}[8]{\draw(#1*1.5-0.375+.15,#2)node{\rotatebox{90}{$\left[\begin{smallmatrix}#3&#4&#5\\#6&#7&#8\end{smallmatrix}\right]$}};}
\nc{\rfrdmat}[8]{\draw(#1*1.5+.375-.15,#2)node{\rotatebox{90}{$\left[\begin{smallmatrix}#3&#4&#5\\#6&#7&#8\end{smallmatrix}\right]$}};}
\nc{\xL}{[x]_{\! _\gL}}\nc{\yL}{[y]_{\! _\gL}}\nc{\xR}{[x]_{\! _\gR}}\nc{\yR}{[y]_{\! _\gR}}\nc{\xH}{[x]_{\! _\gH}}\nc{\yH}{[y]_{\! _\gH}}\nc{\XK}{[X]_{\! _\gK}}\nc{\xK}{[x]_{\! _\gK}}
\nc{\RegSija}{\Reg(\Sija)}
\nc{\MnmK}{\M_{nm}^K}
\nc{\cC}{\mathcal C}
\nc{\cR}{\mathcal R}
\nc{\Ckl}{\cC_k(l)}
\nc{\Rkl}{\cR_k(l)}
\nc{\Cmr}{\cC_m(r)}
\nc{\Rmr}{\cR_m(r)}
\nc{\Cnr}{\cC_n(r)}
\nc{\Rnr}{\cR_n(r)}
\nc{\Z}{\mathbb Z}

\nc{\Reg}{\operatorname{Reg}}
\nc{\RP}{\operatorname{RP}}
\nc{\TXa}{\T_X^a}
\nc{\TXA}{\T(X,A)}
\nc{\TXal}{\T(X,\al)}
\nc{\RegTXa}{\Reg(\TXa)}
\nc{\RegTXA}{\Reg(\TXA)}
\nc{\RegTXal}{\Reg(\TXal)}
\nc{\PalX}{\P_\al(X)}
\nc{\EAX}{\E_A(X)}
\nc{\Bb}{\overline{B}}
\nc{\bb}{\overline{\be}}
\nc{\bw}{{\bf w}}
\nc{\bz}{{\bf z}}
\nc{\TASA}{\T_A\sm\S_A}
\nc{\Ub}{\overline{U}}
\nc{\Vb}{\overline{V}}
\nc{\eb}{\overline{e}}
\nc{\EXa}{\E_X^a}
\nc{\oijr}{1\leq i<j\leq r}
\nc{\veb}{\overline{\ve}}
\nc{\bbT}{\mathbb T}
\nc{\Surj}{\operatorname{Surj}}
\nc{\Sone}{S^{(1)}}
\nc{\fillbox}[2]{\draw[fill=gray!30](#1,#2)--(#1+1,#2)--(#1+1,#2+1)--(#1,#2+1)--(#1,#2);}
\nc{\raa}{\rangle_J}
\nc{\raJ}{\rangle_J}
\nc{\Ea}{E_J}
\nc{\EJ}{E_J}
\nc{\ep}{\epsilon} \nc{\ve}{\varepsilon}
\nc{\IXa}{\I_X^a}
\nc{\RegIXa}{\Reg(\IXa)}
\nc{\JXa}{\J_X^a}
\nc{\RegJXa}{\Reg(\JXa)}
\nc{\IXA}{\I(X,A)}
\nc{\IAX}{\I(A,X)}
\nc{\RegIXA}{\Reg(\IXA)}
\nc{\RegIAX}{\Reg(\IAX)}
\nc{\trans}[2]{\left(\begin{smallmatrix} #1 \\ #2 \end{smallmatrix}\right)}
\nc{\bigtrans}[2]{\left(\begin{matrix} #1 \\ #2 \end{matrix}\right)}
\nc{\lmap}[1]{\mapstochar \xrightarrow {\ #1\ }}
\nc{\EaTXa}{E}

\nc{\gL}{\mathscr L}
\nc{\gR}{\mathscr R}
\nc{\gH}{\mathscr H}
\nc{\gJ}{\mathscr J}
\nc{\gD}{\mathscr D}
\nc{\gK}{\mathscr K}
\nc{\gLa}{\mathscr L^a}
\nc{\gRa}{\mathscr R^a}
\nc{\gHa}{\mathscr H^a}
\nc{\gJa}{\mathscr J^a}
\nc{\gDa}{\mathscr D^a}
\nc{\gKa}{\mathscr K^a}
\nc{\gLJ}{\mathscr L^J}
\nc{\gRJ}{\mathscr R^J}
\nc{\gHJ}{\mathscr H^J}
\nc{\gJJ}{\mathscr J^J}
\nc{\gDJ}{\mathscr D^J}
\nc{\gKJ}{\mathscr K^J}
\nc{\gLh}{\widehat{\mathscr L}^J}
\nc{\gRh}{\widehat{\mathscr R}^J}
\nc{\gHh}{\widehat{\mathscr H}^J}
\nc{\gJh}{\widehat{\mathscr J}^J}
\nc{\gDh}{\widehat{\mathscr D}^J}
\nc{\gKh}{\widehat{\mathscr K}^J}
\nc{\Lh}{\widehat{L}^J}
\nc{\Rh}{\widehat{R}^J}
\nc{\Hh}{\widehat{H}^J}
\nc{\Jh}{\widehat{J}^J}
\nc{\Dh}{\widehat{D}^J}
\nc{\Kh}{\widehat{K}^J}
\nc{\gLb}{\widehat{\mathscr L}}
\nc{\gRb}{\widehat{\mathscr R}}
\nc{\gHb}{\widehat{\mathscr H}}
\nc{\gJb}{\widehat{\mathscr J}}
\nc{\gDb}{\widehat{\mathscr D}}
\nc{\gKb}{\widehat{\mathscr K}}
\nc{\Lb}{\widehat{L}^J}
\nc{\Rb}{\widehat{R}^J}
\nc{\Hb}{\widehat{H}^J}
\nc{\Jb}{\widehat{J}^J}
\nc{\Db}{\overline{D}}
\nc{\Kb}{\widehat{K}}

\hyphenation{mon-oid mon-oids}

\nc{\itemit}[1]{\item[\emph{(#1)}]}
\nc{\E}{\mathcal E}
\nc{\TX}{\T(X)}
\nc{\TXP}{\T(X,\P)}
\nc{\EX}{\E(X)}
\nc{\EXP}{\E(X,\P)}
\nc{\SX}{\S(X)}
\nc{\SXP}{\S(X,\P)}
\nc{\Sing}{\operatorname{Sing}}
\nc{\idrank}{\operatorname{idrank}}
\nc{\SingXP}{\Sing(X,\P)}
\nc{\De}{\Delta}
\nc{\sgp}{\operatorname{sgp}}
\nc{\mon}{\operatorname{mon}}
\nc{\Dn}{\mathcal D_n}
\nc{\Dm}{\mathcal D_m}

\nc{\lline}[1]{\draw(3*#1,0)--(3*#1+2,0);}
\nc{\uline}[1]{\draw(3*#1,5)--(3*#1+2,5);}
\nc{\thickline}[2]{\draw(3*#1,5)--(3*#2,0); \draw(3*#1+2,5)--(3*#2+2,0) ;}
\nc{\thicklabel}[3]{\draw(3*#1+1+3*#2*0.15-3*#1*0.15,4.25)node{{\tiny $#3$}};}

\nc{\slline}[3]{\draw(3*#1+#3,0+#2)--(3*#1+2+#3,0+#2);}
\nc{\suline}[3]{\draw(3*#1+#3,5+#2)--(3*#1+2+#3,5+#2);}
\nc{\sthickline}[4]{\draw(3*#1+#4,5+#3)--(3*#2+#4,0+#3); \draw(3*#1+2+#4,5+#3)--(3*#2+2+#4,0+#3) ;}
\nc{\sthicklabel}[5]{\draw(3*#1+1+3*#2*0.15-3*#1*0.15+#5,4.25+#4)node{{\tiny $#3$}};}

\nc{\stll}[5]{\sthickline{#1}{#2}{#4}{#5} \sthicklabel{#1}{#2}{#3}{#4}{#5}}
\nc{\tll}[3]{\stll{#1}{#2}{#3}00}

\nc{\mfourpic}[9]{
\slline1{#9}0
\slline3{#9}0
\slline4{#9}0
\slline5{#9}0
\suline1{#9}0
\suline3{#9}0
\suline4{#9}0
\suline5{#9}0
\stll1{#1}{#5}{#9}{0}
\stll3{#2}{#6}{#9}{0}
\stll4{#3}{#7}{#9}{0}
\stll5{#4}{#8}{#9}{0}
\draw[dotted](6,0+#9)--(8,0+#9);
\draw[dotted](6,5+#9)--(8,5+#9);
}
\nc{\vdotted}[1]{
\draw[dotted](3*#1,10)--(3*#1,15);
\draw[dotted](3*#1+2,10)--(3*#1+2,15);
}

\nc{\Clab}[2]{
\sthicklabel{#1}{#1}{{}_{\phantom{#1}}C_{#1}}{1.25+5*#2}0
}
\nc{\sClab}[3]{
\sthicklabel{#1}{#1}{{}_{\phantom{#1}}C_{#1}}{1.25+5*#2}{#3}
}
\nc{\Clabl}[3]{
\sthicklabel{#1}{#1}{{}_{\phantom{#3}}C_{#3}}{1.25+5*#2}0
}
\nc{\sClabl}[4]{
\sthicklabel{#1}{#1}{{}_{\phantom{#4}}C_{#4}}{1.25+5*#2}{#3}
}
\nc{\Clabll}[3]{
\sthicklabel{#1}{#1}{C_{#3}}{1.25+5*#2}0
}
\nc{\sClabll}[4]{
\sthicklabel{#1}{#1}{C_{#3}}{1.25+5*#2}{#3}
}

\nc{\mtwopic}[6]{
\slline1{#6*5}{#5}
\slline2{#6*5}{#5}
\suline1{#6*5}{#5}
\suline2{#6*5}{#5}
\stll1{#1}{#3}{#6*5}{#5}
\stll2{#2}{#4}{#6*5}{#5}
}
\nc{\mtwopicl}[6]{
\slline1{#6*5}{#5}
\slline2{#6*5}{#5}
\suline1{#6*5}{#5}
\suline2{#6*5}{#5}
\stll1{#1}{#3}{#6*5}{#5}
\stll2{#2}{#4}{#6*5}{#5}
\sClabl1{#6}{#5}{i}
\sClabl2{#6}{#5}{j}
}

\nc{\keru}{\operatorname{ker}^\wedge} \nc{\kerl}{\operatorname{ker}_\vee}

\nc{\coker}{\operatorname{coker}}
\nc{\KER}{\ker}
\nc{\N}{\mathbb N}
\nc{\LaBn}{L_\al(\B_n)}
\nc{\RaBn}{R_\al(\B_n)}
\nc{\LaPBn}{L_\al(\PB_n)}
\nc{\RaPBn}{R_\al(\PB_n)}
\nc{\rhorBn}{\rho_r(\B_n)}
\nc{\DrBn}{D_r(\B_n)}
\nc{\DrPn}{D_r(\P_n)}
\nc{\DrPBn}{D_r(\PB_n)}
\nc{\DrKn}{D_r(\K_n)}
\nc{\alb}{\al_{\vee}}
\nc{\beb}{\be^{\wedge}}
\nc{\Bal}{\operatorname{Bal}}
\nc{\Red}{\operatorname{Red}}
\nc{\Pnxi}{\P_n^\xi}
\nc{\Bnxi}{\B_n^\xi}
\nc{\PBnxi}{\PB_n^\xi}
\nc{\Knxi}{\K_n^\xi}
\nc{\C}{\mathscr C}
\nc{\exi}{e^\xi}
\nc{\Exi}{E^\xi}
\nc{\eximu}{e^\xi_\mu}
\nc{\Eximu}{E^\xi_\mu}
\nc{\REF}{ {\red [Ref?]} }
\nc{\GL}{\operatorname{GL}}
\rnc{\O}{\mathcal O}

\nc{\vtx}[2]{\fill (#1,#2)circle(.2);}
\nc{\lvtx}[2]{\fill (#1,0)circle(.2);}
\nc{\uvtx}[2]{\fill (#1,1.5)circle(.2);}

\nc{\Eq}{\mathfrak{Eq}}
\nc{\Gau}{\Ga^\wedge} \nc{\Gal}{\Ga_\vee}
\nc{\Lamu}{\Lam^\wedge} \nc{\Laml}{\Lam_\vee}
\nc{\bX}{{\bf X}}
\nc{\bY}{{\bf Y}}
\nc{\ds}{\displaystyle}

\nc{\uuvert}[1]{\fill (#1,3)circle(.2);}
\nc{\uuuvert}[1]{\fill (#1,4.5)circle(.2);}
\nc{\overt}[1]{\fill (#1,0)circle(.1);}
\nc{\overtl}[3]{\node[vertex] (#3) at (#1,0) {  {\tiny $#2$} };}
\nc{\cv}[2]{\draw(#1,1.5) to [out=270,in=90] (#2,0);}
\nc{\cvs}[2]{\draw(#1,1.5) to [out=270+30,in=90+30] (#2,0);}
\nc{\ucv}[2]{\draw(#1,3) to [out=270,in=90] (#2,1.5);}
\nc{\uucv}[2]{\draw(#1,4.5) to [out=270,in=90] (#2,3);}
\nc{\textpartn}[1]{{\lower1.0 ex\hbox{\begin{tikzpicture}[xscale=.3,yscale=0.3] #1 \end{tikzpicture}}}}
\nc{\textpartnx}[2]{{\lower1.0 ex\hbox{\begin{tikzpicture}[xscale=.3,yscale=0.3] 
\foreach \x in {1,...,#1}
{ \uvert{\x} \lvert{\x} }
#2 \end{tikzpicture}}}}
\nc{\disppartnx}[2]{{\lower1.0 ex\hbox{\begin{tikzpicture}[scale=0.3] 
\foreach \x in {1,...,#1}
{ \uvert{\x} \lvert{\x} }
#2 \end{tikzpicture}}}}
\nc{\disppartnxd}[2]{{\lower2.1 ex\hbox{\begin{tikzpicture}[scale=0.3] 
\foreach \x in {1,...,#1}
{ \uuvert{\x} \uvert{\x} \lvert{\x} }
#2 \end{tikzpicture}}}}
\nc{\disppartnxdn}[2]{{\lower2.1 ex\hbox{\begin{tikzpicture}[scale=0.3] 
\foreach \x in {1,...,#1}
{ \uuvert{\x} \lvert{\x} }
#2 \end{tikzpicture}}}}
\nc{\disppartnxdd}[2]{{\lower3.6 ex\hbox{\begin{tikzpicture}[scale=0.3] 
\foreach \x in {1,...,#1}
{ \uuuvert{\x} \uuvert{\x} \uvert{\x} \lvert{\x} }
#2 \end{tikzpicture}}}}

\nc{\dispgax}[2]{{\lower0.0 ex\hbox{\begin{tikzpicture}[scale=0.3] 
#2
\foreach \x in {1,...,#1}
{\lvert{\x} }
 \end{tikzpicture}}}}
\nc{\textgax}[2]{{\lower0.4 ex\hbox{\begin{tikzpicture}[scale=0.3] 
#2
\foreach \x in {1,...,#1}
{\lvert{\x} }
 \end{tikzpicture}}}}
\nc{\textlinegraph}[2]{{\raise#1 ex\hbox{\begin{tikzpicture}[scale=0.8] 
#2
 \end{tikzpicture}}}}
\nc{\textlinegraphl}[2]{{\raise#1 ex\hbox{\begin{tikzpicture}[scale=0.8] 
\tikzstyle{vertex}=[circle,draw=black, fill=white, inner sep = 0.07cm]
#2
 \end{tikzpicture}}}}
\nc{\displinegraph}[1]{{\lower0.0 ex\hbox{\begin{tikzpicture}[scale=0.6] 
#1
 \end{tikzpicture}}}}
 
\nc{\disppartnthreeone}[1]{{\lower1.0 ex\hbox{\begin{tikzpicture}[scale=0.3] 
\foreach \x in {1,2,3,5,6}
{ \uvert{\x} }
\foreach \x in {1,2,4,5,6}
{ \lvert{\x} }
\draw[dotted] (3.5,1.5)--(4.5,1.5);
\draw[dotted] (2.5,0)--(3.5,0);
#1 \end{tikzpicture}}}}

\nc{\partn}[4]{\left( \begin{array}{c|c} 
#1 \ & \ #3 \ \ \\ \cline{2-2}
#2 \ & \ #4 \ \
\end{array} \!\!\! \right)}
\nc{\partnlong}[6]{\partn{#1}{#2}{#3,\ #4}{#5,\ #6}} 
\nc{\partnsh}[2]{\left( \begin{array}{c} 
#1 \\
#2 
\end{array} \right)}
\nc{\partncodefz}[3]{\partn{#1}{#2}{#3}{\emptyset}}
\nc{\partndefz}[3]{{\partn{#1}{#2}{\emptyset}{#3}}}
\nc{\partnlast}[2]{\left( \begin{array}{c|c}
#1 \ &  \ #2 \\
#1 \ &  \ #2
\end{array} \right)}

\nc{\uuarcx}[3]{\draw(#1,3)arc(180:270:#3) (#1+#3,3-#3)--(#2-#3,3-#3) (#2-#3,3-#3) arc(270:360:#3);}
\nc{\uuarc}[2]{\uuarcx{#1}{#2}{.4}}
\nc{\uuuarcx}[3]{\draw(#1,4.5)arc(180:270:#3) (#1+#3,4.5-#3)--(#2-#3,4.5-#3) (#2-#3,4.5-#3) arc(270:360:#3);}
\nc{\uuuarc}[2]{\uuuarcx{#1}{#2}{.4}}
\nc{\udarcx}[3]{\draw(#1,1.5)arc(180:90:#3) (#1+#3,1.5+#3)--(#2-#3,1.5+#3) (#2-#3,1.5+#3) arc(90:0:#3);}
\nc{\udarc}[2]{\udarcx{#1}{#2}{.4}}
\nc{\uudarcx}[3]{\draw(#1,3)arc(180:90:#3) (#1+#3,3+#3)--(#2-#3,3+#3) (#2-#3,3+#3) arc(90:0:#3);}
\nc{\uudarc}[2]{\uudarcx{#1}{#2}{.4}}
\nc{\darcxhalf}[3]{\draw(#1,0)arc(180:90:#3) (#1+#3,#3)--(#2,#3) ;}
\nc{\darchalf}[2]{\darcxhalf{#1}{#2}{.4}}
\nc{\uarcxhalf}[3]{\draw(#1,2)arc(180:270:#3) (#1+#3,2-#3)--(#2,2-#3) ;}
\nc{\uarchalf}[2]{\uarcxhalf{#1}{#2}{.4}}
\nc{\uarcxhalfr}[3]{\draw (#1,2-#3)--(#2-#3,2-#3) (#2-#3,2-#3) arc(270:360:#3);}
\nc{\uarchalfr}[2]{\uarcxhalfr{#1}{#2}{.4}}
\nc{\darcxhalfr}[3]{\draw (#1,#3)--(#2-#3,#3) (#2-#3,#3) arc(90:0:#3);}
\nc{\darchalfr}[2]{\darcxhalfr{#1}{#2}{.4}}

\nc{\bdarcx}[3]{\draw[blue](#1,0)arc(180:90:#3) (#1+#3,#3)--(#2-#3,#3) (#2-#3,#3) arc(90:0:#3);}
\nc{\bdarc}[2]{\darcx{#1}{#2}{.4}}
\nc{\rduarcx}[3]{\draw[red](#1,0)arc(180:270:#3) (#1+#3,0-#3)--(#2-#3,0-#3) (#2-#3,0-#3) arc(270:360:#3);}
\nc{\rduarc}[2]{\uarcx{#1}{#2}{.4}}
\nc{\duarcx}[3]{\draw(#1,0)arc(180:270:#3) (#1+#3,0-#3)--(#2-#3,0-#3) (#2-#3,0-#3) arc(270:360:#3);}
\nc{\duarc}[2]{\uarcx{#1}{#2}{.4}}

\nc{\uuv}[1]{\fill (#1,4)circle(.1);}
\nc{\uv}[1]{\fill (#1,2)circle(.1);}
\nc{\lv}[1]{\fill (#1,0)circle(.1);}
\nc{\uvred}[1]{\fill[red] (#1,2)circle(.1);}
\nc{\lvred}[1]{\fill[red] (#1,0)circle(.1);}
\nc{\lvwhite}[1]{\fill[white] (#1,0)circle(.1);}

\nc{\uvs}[1]{{
\foreach \x in {#1}
{ \uv{\x}}
}}
\nc{\uuvs}[1]{{
\foreach \x in {#1}
{ \uuv{\x}}
}}
\nc{\lvs}[1]{{
\foreach \x in {#1}
{ \lv{\x}}
}}

\nc{\uvreds}[1]{{
\foreach \x in {#1}
{ \uvred{\x}}
}}
\nc{\lvreds}[1]{{
\foreach \x in {#1}
{ \lvred{\x}}
}}

\nc{\uudotted}[2]{\draw [dotted] (#1,4)--(#2,4);}
\nc{\uudotteds}[1]{{
\foreach \x/\y in {#1}
{ \uudotted{\x}{\y}}
}}
\nc{\uudottedsm}[2]{\draw [dotted] (#1+.4,4)--(#2-.4,4);}
\nc{\uudottedsms}[1]{{
\foreach \x/\y in {#1}
{ \uudottedsm{\x}{\y}}
}}
\nc{\udottedsm}[2]{\draw [dotted] (#1+.4,2)--(#2-.4,2);}
\nc{\udottedsms}[1]{{
\foreach \x/\y in {#1}
{ \udottedsm{\x}{\y}}
}}
\nc{\udotted}[2]{\draw [dotted] (#1+.5,2)--(#2-.5,2);}
\nc{\udotteds}[1]{{
\foreach \x/\y in {#1}
{ \udotted{\x}{\y}}
}}
\nc{\ldotted}[2]{\draw [dotted] (#1+.5,0)--(#2-.5,0);}
\nc{\ldotteds}[1]{{
\foreach \x/\y in {#1}
{ \ldotted{\x}{\y}}
}}
\nc{\ldottedsm}[2]{\draw [dotted] (#1+.4,0)--(#2-.4,0);}
\nc{\ldottedsms}[1]{{
\foreach \x/\y in {#1}
{ \ldottedsm{\x}{\y}}
}}

\nc{\stlinest}[2]{\draw(#1,4)--(#2,0);}

\nc{\stlined}[2]{\draw[dotted](#1,2)--(#2,0);}

\nc{\tlab}[2]{\draw(#1,2)node[above]{\tiny $#2$};}
\nc{\tudots}[1]{\draw(#1,2)node{$\cdots$};}
\nc{\tldots}[1]{\draw(#1,0)node{$\cdots$};}

\nc{\uvw}[1]{\fill[white] (#1,2)circle(.1);}
\nc{\huv}[1]{\fill (#1,1)circle(.1);}
\nc{\llv}[1]{\fill (#1,-2)circle(.1);}
\nc{\arcup}[2]{
\draw(#1,2)arc(180:270:.4) (#1+.4,1.6)--(#2-.4,1.6) (#2-.4,1.6) arc(270:360:.4);
}
\nc{\harcup}[2]{
\draw(#1,1)arc(180:270:.4) (#1+.4,.6)--(#2-.4,.6) (#2-.4,.6) arc(270:360:.4);
}
\nc{\arcdn}[2]{
\draw(#1,0)arc(180:90:.4) (#1+.4,.4)--(#2-.4,.4) (#2-.4,.4) arc(90:0:.4);
}
\nc{\cve}[2]{
\draw(#1,2) to [out=270,in=0] (0.5*#1+0.5*#2,1) to [out=180,in=90] (#2,0);
}
\nc{\hcve}[2]{
\draw(#1,1) to [out=270,in=90] (#2,0);
}
\nc{\catarc}[3]{
\draw(#1,2)arc(180:270:#3) (#1+#3,2-#3)--(#2-#3,2-#3) (#2-#3,2-#3) arc(270:360:#3);
}

\nc{\arcr}[2]{
\draw[red](#1,0)arc(180:90:.4) (#1+.4,.4)--(#2-.4,.4) (#2-.4,.4) arc(90:0:.4);
}
\nc{\arcb}[2]{
\draw[blue](#1,2-2)arc(180:270:.4) (#1+.4,1.6-2)--(#2-.4,1.6-2) (#2-.4,1.6-2) arc(270:360:.4);
}
\nc{\loopr}[1]{
\draw[blue](#1,-2) edge [out=130,in=50,loop] ();
}
\nc{\loopb}[1]{
\draw[red](#1,-2) edge [out=180+130,in=180+50,loop] ();
}
\nc{\redto}[2]{\draw[red](#1,0)--(#2,0);}
\nc{\bluto}[2]{\draw[blue](#1,0)--(#2,0);}
\nc{\dotto}[2]{\draw[dotted](#1,0)--(#2,0);}
\nc{\lloopr}[2]{\draw[red](#1,0)arc(0:360:#2);}
\nc{\lloopb}[2]{\draw[blue](#1,0)arc(0:360:#2);}
\nc{\rloopr}[2]{\draw[red](#1,0)arc(-180:180:#2);}
\nc{\rloopb}[2]{\draw[blue](#1,0)arc(-180:180:#2);}
\nc{\uloopr}[2]{\draw[red](#1,0)arc(-270:270:#2);}
\nc{\uloopb}[2]{\draw[blue](#1,0)arc(-270:270:#2);}
\nc{\dloopr}[2]{\draw[red](#1,0)arc(-90:270:#2);}
\nc{\dloopb}[2]{\draw[blue](#1,0)arc(-90:270:#2);}
\nc{\llloopr}[2]{\draw[red](#1,0-2)arc(0:360:#2);}
\nc{\llloopb}[2]{\draw[blue](#1,0-2)arc(0:360:#2);}
\nc{\lrloopr}[2]{\draw[red](#1,0-2)arc(-180:180:#2);}
\nc{\lrloopb}[2]{\draw[blue](#1,0-2)arc(-180:180:#2);}
\nc{\ldloopr}[2]{\draw[red](#1,0-2)arc(-270:270:#2);}
\nc{\ldloopb}[2]{\draw[blue](#1,0-2)arc(-270:270:#2);}
\nc{\luloopr}[2]{\draw[red](#1,0-2)arc(-90:270:#2);}
\nc{\luloopb}[2]{\draw[blue](#1,0-2)arc(-90:270:#2);}

\nc{\larcb}[2]{
\draw[blue](#1,0-2)arc(180:90:.4) (#1+.4,.4-2)--(#2-.4,.4-2) (#2-.4,.4-2) arc(90:0:.4);
}
\nc{\larcr}[2]{
\draw[red](#1,2-2-2)arc(180:270:.4) (#1+.4,1.6-2-2)--(#2-.4,1.6-2-2) (#2-.4,1.6-2-2) arc(270:360:.4);
}

\rnc{\H}{\mathrel{\mathscr H}}
\rnc{\L}{\mathrel{\mathscr L}}
\nc{\R}{\mathrel{\mathscr R}}
\nc{\D}{\mathrel{\mathscr D}}
\nc{\J}{\mathrel{\mathscr J}}
\nc{\Ht}{\mathrel{\mathscr H^\tau}}
\nc{\Lt}{\mathrel{\mathscr L^\tau}}
\nc{\Rt}{\mathrel{\mathscr R^\tau}}
\nc{\Dt}{\mathrel{\mathscr D^\tau}}
\nc{\Jt}{\mathrel{\mathscr J^\tau}}
\nc{\Kt}{\mathrel{\mathscr K^\tau}}

\nc{\leqR}{\leq_{\R}}
\nc{\geqR}{\geq_{\R}}
\nc{\leqL}{\leq_{\L}}
\nc{\leqJ}{\leq_{\J}}
\nc{\leqK}{\leq_{\K}}
\nc{\leqRt}{\leq_{\R}^\tau}
\nc{\geqRt}{\geq_{\R}^\tau}
\nc{\geqJt}{\geq_{\J}^\tau}
\nc{\leqLt}{\leq_{\L}^\tau}
\nc{\leqJt}{\leq_{\J}^\tau}
\nc{\leqKt}{\leq_{\K}^\tau}

\nc{\ssim}{\mathrel{\raise0.25 ex\hbox{\oalign{$\approx$\crcr\noalign{\kern-0.84 ex}$\approx$}}}}
\nc{\POI}{\mathcal{O}}
\nc{\wb}{\overline{w}}
\nc{\ub}{\overline{u}}
\nc{\vb}{\overline{v}}
\nc{\fb}{\overline{f}}
\nc{\gb}{\overline{g}}
\nc{\hb}{\overline{h}}
\nc{\pb}{\overline{p}}
\rnc{\sb}{\overline{s}}
\nc{\XR}{\pres{X}{R\,}}
\nc{\YQ}{\pres{Y}{Q}}
\nc{\ZP}{\pres{Z}{P\,}}
\nc{\XRone}{\pres{X_1}{R_1}}
\nc{\XRtwo}{\pres{X_2}{R_2}}
\nc{\XRthree}{\pres{X_1\cup X_2}{R_1\cup R_2\cup R_3}}
\nc{\er}{\eqref}
\nc{\larr}{\mathrel{\hspace{-0.35 ex}>\hspace{-1.1ex}-}\hspace{-0.35 ex}}
\nc{\rarr}{\mathrel{\hspace{-0.35 ex}-\hspace{-0.5ex}-\hspace{-2.3ex}>\hspace{-0.35 ex}}}
\nc{\lrarr}{\mathrel{\hspace{-0.35 ex}>\hspace{-1.1ex}-\hspace{-0.5ex}-\hspace{-2.3ex}>\hspace{-0.35 ex}}}
\nc{\nn}{\tag*{}}
\nc{\epfal}{\tag*{$\Box$}}
\nc{\tagd}[1]{\tag*{(#1)$'$}}
\nc{\ldb}{[\![}
\nc{\rdb}{]\!]}
\nc{\sm}{\setminus}
\nc{\I}{\mathcal I}
\nc{\InSn}{\I_n\setminus\S_n}
\nc{\dom}{\operatorname{dom}} \nc{\codom}{\operatorname{codom}}
\nc{\ojin}{1\leq j<i\leq n}
\nc{\eh}{\widehat{e}}
\nc{\wh}{\widehat{w}}
\nc{\uh}{\widehat{u}}
\nc{\vh}{\widehat{v}}
\nc{\sh}{\widehat{s}}
\nc{\fh}{\widehat{f}}
\nc{\textres}[1]{\text{\emph{#1}}}
\nc{\aand}{\emph{\ and \ }}
\nc{\iif}{\emph{\ if \ }}
\nc{\textlarr}{\ \larr\ }
\nc{\textrarr}{\ \rarr\ }
\nc{\textlrarr}{\ \lrarr\ }

\nc{\comma}{,\ }

\nc{\COMMA}{,\quad}
\nc{\TnSn}{\T_n\setminus\S_n} 
\nc{\TmSm}{\T_m\setminus\S_m} 
\nc{\TXSX}{\T_X\setminus\S_X} 
\rnc{\S}{\mathcal S}

\nc{\T}{\mathcal T} 
\nc{\A}{\mathscr A} 
\nc{\B}{\mathcal B} 
\rnc{\P}{\mathcal P} 
\nc{\K}{\mathrel{\mathscr K}}
\nc{\PB}{\mathcal{PB}} 
\nc{\rank}{\operatorname{rank}}

\nc{\mtt}{\!\!\!\mt\!\!\!}

\nc{\sub}{\subseteq}
\nc{\la}{\langle}
\nc{\ra}{\rangle}
\nc{\mt}{\mapsto}
\nc{\im}{\mathrm{im}}
\nc{\id}{\mathrm{id}}
\nc{\al}{\alpha}
\nc{\be}{\beta}
\nc{\ga}{\gamma}
\nc{\Ga}{\Gamma}
\nc{\de}{\delta}
\nc{\ka}{\kappa}
\nc{\lam}{\lambda}
\nc{\Lam}{\Lambda}
\nc{\si}{\sigma}
\nc{\Si}{\Sigma}
\nc{\oijn}{1\leq i<j\leq n}
\nc{\oijm}{1\leq i<j\leq m}

\nc{\comm}{\rightleftharpoons}
\nc{\AND}{\qquad\text{and}\qquad}

\nc{\bit}{\vspace{-3 truemm}\begin{itemize}}
\nc{\bitbmc}{\begin{itemize}\begin{multicols}}
\nc{\bmc}{\begin{itemize}\begin{multicols}}
\nc{\emc}{\end{multicols}\end{itemize}\vspace{-3 truemm}}
\nc{\eit}{\end{itemize}\vspace{-3 truemm}}
\nc{\ben}{\vspace{-3 truemm}\begin{enumerate}}
\nc{\een}{\end{enumerate}\vspace{-3 truemm}}
\nc{\eitres}{\end{itemize}}

\nc{\set}[2]{\{ {#1} : {#2} \}} 
\nc{\bigset}[2]{\big\{ {#1}: {#2} \big\}} 
\nc{\Bigset}[2]{\left\{ \,{#1} :{#2}\, \right\}}

\nc{\pres}[2]{\la {#1} \,|\, {#2} \ra}
\nc{\bigpres}[2]{\big\la {#1} \,\big|\, {#2} \big\ra}
\nc{\Bigpres}[2]{\Big\la \,{#1}\, \,\Big|\, \,{#2}\, \Big\ra}
\nc{\Biggpres}[2]{\Bigg\la {#1} \,\Bigg|\, {#2} \Bigg\ra}

\nc{\pf}{\noindent{\bf Proof.}  }
\nc{\epf}{\hfill$\Box$\bigskip}
\nc{\epfres}{\hfill$\Box$}
\nc{\pfnb}{\pf}
\nc{\epfnb}{\bigskip}
\nc{\pfthm}[1]{\bigskip \noindent{\bf Proof of Theorem \ref{#1}}\,\,  } 
\nc{\pfprop}[1]{\bigskip \noindent{\bf Proof of Proposition \ref{#1}}\,\,  } 
\nc{\epfreseq}{\tag*{$\Box$}}

\nc{\uvert}[1]{\fill (#1,2)circle(.2);}
\rnc{\lvert}[1]{\fill (#1,0)circle(.2);}
\nc{\uvertcol}[2]{\fill[#2] (#1,2)circle(.2);}
\nc{\lvertcol}[2]{\fill[#2] (#1,0)circle(.2);}
\nc{\guvert}[1]{\fill[lightgray] (#1,2)circle(.2);}
\nc{\glvert}[1]{\fill[lightgray] (#1,0)circle(.2);}
\nc{\uvertx}[2]{\fill (#1,#2)circle(.2);}
\nc{\guvertx}[2]{\fill[lightgray] (#1,#2)circle(.2);}
\nc{\uvertxs}[2]{
\foreach \x in {#1}
{ \uvertx{\x}{#2}  }
}
\nc{\guvertxs}[2]{
\foreach \x in {#1}
{ \guvertx{\x}{#2}  }
}

\nc{\uvertth}[2]{\fill (#1,2)circle(#2);}
\nc{\lvertth}[2]{\fill (#1,0)circle(#2);}
\nc{\uvertths}[2]{
\foreach \x in {#1}
{ \uvertth{\x}{#2}  }
}
\nc{\lvertths}[2]{
\foreach \x in {#1}
{ \lvertth{\x}{#2}  }
}

\nc{\vertlabel}[2]{\draw(#1,2+.5)node{{\tiny $#2$}};}
\nc{\vertlabelh}[2]{\draw(#1,2+.4)node{{\tiny $#2$}};}
\nc{\vertlabelhh}[2]{\draw(#1,2+.6)node{{\tiny $#2$}};}
\nc{\vertlabelhhh}[2]{\draw(#1,2+.64)node{{\tiny $#2$}};}
\nc{\vertlabelup}[2]{\draw(#1,4+.6)node{{\tiny $#2$}};}
\nc{\vertlabels}[1]{
{\foreach \x/\y in {#1}
{ \vertlabel{\x}{\y} }
}
}

\nc{\dvertlabel}[2]{\draw(#1,-.4)node{{\tiny $#2$}};}
\nc{\dvertlabels}[1]{
{\foreach \x/\y in {#1}
{ \dvertlabel{\x}{\y} }
}
}
\nc{\vertlabelsh}[1]{
{\foreach \x/\y in {#1}
{ \vertlabelh{\x}{\y} }
}
}
\nc{\vertlabelshh}[1]{
{\foreach \x/\y in {#1}
{ \vertlabelhh{\x}{\y} }
}
}
\nc{\vertlabelshhh}[1]{
{\foreach \x/\y in {#1}
{ \vertlabelhhh{\x}{\y} }
}
}

\nc{\vertlabelx}[3]{\draw(#1,2+#3+.6)node{{\tiny $#2$}};}
\nc{\vertlabelxs}[2]{
{\foreach \x/\y in {#1}
{ \vertlabelx{\x}{\y}{#2} }
}
}

\nc{\vertlabelupdash}[2]{\draw(#1,2.7)node{{\tiny $\phantom{'}#2'$}};}
\nc{\vertlabelupdashess}[1]{
{\foreach \x/\y in {#1}
{\vertlabelupdash{\x}{\y}}
}
}

\nc{\vertlabeldn}[2]{\draw(#1,0-.6)node{{\tiny $\phantom{'}#2'$}};}
\nc{\vertlabeldnph}[2]{\draw(#1,0-.6)node{{\tiny $\phantom{'#2'}$}};}

\nc{\vertlabelups}[1]{
{\foreach \x in {#1}
{\vertlabel{\x}{\x}}
}
}
\nc{\vertlabeldns}[1]{
{\foreach \x in {#1}
{\vertlabeldn{\x}{\x}}
}
}
\nc{\vertlabeldnsph}[1]{
{\foreach \x in {#1}
{\vertlabeldnph{\x}{\x}}
}
}

\nc{\dotsup}[2]{\draw [dotted] (#1+.6,2)--(#2-.6,2);}
\nc{\dotsupx}[3]{\draw [dotted] (#1+.6,#3)--(#2-.6,#3);}
\nc{\dotsdn}[2]{\draw [dotted] (#1+.6,0)--(#2-.6,0);}
\nc{\dotsups}[1]{\foreach \x/\y in {#1}
{ \dotsup{\x}{\y} }
}
\nc{\dotsupxs}[2]{\foreach \x/\y in {#1}
{ \dotsupx{\x}{\y}{#2} }
}
\nc{\dotsdns}[1]{\foreach \x/\y in {#1}
{ \dotsdn{\x}{\y} }
}

\nc{\nodropcustpartn}[3]{
\begin{tikzpicture}[scale=.3]
\foreach \x in {#1}
{ \uvert{\x}  }
\foreach \x in {#2}
{ \lvert{\x}  }
#3 \end{tikzpicture}
}

\nc{\custpartn}[3]{{\lower1.4 ex\hbox{
\begin{tikzpicture}[scale=.3]
\foreach \x in {#1}
{ \uvert{\x}  }
\foreach \x in {#2}
{ \lvert{\x}  }
#3 \end{tikzpicture}
}}}

\nc{\smcustpartn}[3]{{\lower0.7 ex\hbox{
\begin{tikzpicture}[scale=.2]
\foreach \x in {#1}
{ \uvert{\x}  }
\foreach \x in {#2}
{ \lvert{\x}  }
#3 \end{tikzpicture}
}}}

\nc{\dropcustpartn}[3]{{\lower5.2 ex\hbox{
\begin{tikzpicture}[scale=.3]
\foreach \x in {#1}
{ \uvert{\x}  }
\foreach \x in {#2}
{ \lvert{\x}  }
#3 \end{tikzpicture}
}}}

\nc{\dropcustpartnx}[4]{{\lower#4 ex\hbox{
\begin{tikzpicture}[scale=.4]
\foreach \x in {#1}
{ \uvert{\x}  }
\foreach \x in {#2}
{ \lvert{\x}  }
#3 \end{tikzpicture}
}}}

\nc{\dropcustpartnxy}[3]{\dropcustpartnx{#1}{#2}{#3}{4.6}}

\nc{\uvertsm}[1]{\fill (#1,2)circle(.15);}
\nc{\lvertsm}[1]{\fill (#1,0)circle(.15);}
\nc{\vertsm}[2]{\fill (#1,#2)circle(.15);}

\nc{\bigdropcustpartn}[3]{{\lower6.93 ex\hbox{
\begin{tikzpicture}[scale=.6]
\foreach \x in {#1}
{ \uvertsm{\x}  }
\foreach \x in {#2}
{ \lvertsm{\x}  }
#3 \end{tikzpicture}
}}}

\nc{\gcustpartn}[5]{{\lower1.4 ex\hbox{
\begin{tikzpicture}[scale=.3]
\foreach \x in {#1}
{ \uvert{\x}  }
\foreach \x in {#2}
{ \guvert{\x}  }
\foreach \x in {#3}
{ \lvert{\x}  }
\foreach \x in {#4}
{ \glvert{\x}  }
#5 \end{tikzpicture}
}}}

\nc{\gcustpartndash}[5]{{\lower6.97 ex\hbox{
\begin{tikzpicture}[scale=.3]
\foreach \x in {#1}
{ \uvert{\x}  }
\foreach \x in {#2}
{ \guvert{\x}  }
\foreach \x in {#3}
{ \lvert{\x}  }
\foreach \x in {#4}
{ \glvert{\x}  }
#5 \end{tikzpicture}
}}}

\nc{\stline}[2]{\draw(#1,2)--(#2,0);}
\nc{\stlines}[1]{
{\foreach \x/\y in {#1}
{ \stline{\x}{\y} }
}
}

\nc{\uarcs}[1]{
{\foreach \x/\y in {#1}
{ \uarc{\x}{\y} }
}
}
\nc{\darcs}[1]{
{\foreach \x/\y in {#1}
{ \darc{\x}{\y} }
}
}

\nc{\stlinests}[1]{
{\foreach \x/\y in {#1}
{ \stlinest{\x}{\y} }
}
}

\nc{\stlineds}[1]{
{\foreach \x/\y in {#1}
{ \stlined{\x}{\y} }
}
}

\nc{\gstline}[2]{\draw[lightgray](#1,2)--(#2,0);}
\nc{\gstlines}[1]{
{\foreach \x/\y in {#1}
{ \gstline{\x}{\y} }
}
}

\nc{\gstlinex}[3]{\draw[lightgray](#1,2+#3)--(#2,0+#3);}
\nc{\gstlinexs}[2]{
{\foreach \x/\y in {#1}
{ \gstlinex{\x}{\y}{#2} }
}
}

\nc{\stlinex}[3]{\draw(#1,2+#3)--(#2,0+#3);}
\nc{\stlinexs}[2]{
{\foreach \x/\y in {#1}
{ \stlinex{\x}{\y}{#2} }
}
}

\nc{\darcx}[3]{\draw(#1,0)arc(180:90:#3) (#1+#3,#3)--(#2-#3,#3) (#2-#3,#3) arc(90:0:#3);}
\nc{\darc}[2]{\darcx{#1}{#2}{.4}}
\nc{\uarcx}[3]{\draw(#1,2)arc(180:270:#3) (#1+#3,2-#3)--(#2-#3,2-#3) (#2-#3,2-#3) arc(270:360:#3);}
\nc{\uarc}[2]{\uarcx{#1}{#2}{.4}}

\nc{\darcxx}[4]{\draw(#1,0+#4)arc(180:90:#3) (#1+#3,#3+#4)--(#2-#3,#3+#4) (#2-#3,#3+#4) arc(90:0:#3);}
\nc{\uarcxx}[4]{\draw(#1,2+#4)arc(180:270:#3) (#1+#3,2-#3+#4)--(#2-#3,2-#3+#4) (#2-#3,2-#3+#4) arc(270:360:#3);}

\makeatletter
\newcommand\footnoteref[1]{\protected@xdef\@thefnmark{\ref{#1}}\@footnotemark}
\makeatother

\newcounter{theorem}
\numberwithin{theorem}{section}

\newtheorem{thm}[theorem]{Theorem}
\newtheorem{lemma}[theorem]{Lemma}
\newtheorem{cor}[theorem]{Corollary}
\newtheorem{prop}[theorem]{Proposition}

\theoremstyle{definition}

\newtheorem{rem}[theorem]{Remark}
\newtheorem{defn}[theorem]{Definition}
\newtheorem{eg}[theorem]{Example}
\newtheorem{ass}[theorem]{Assumption}

\title{Twisted Brauer monoids}

\date{}

\author{
Igor Dolinka\footnote{Department of Mathematics and Informatics, University of Novi Sad, Trg Dositeja Obradovi\'ca 4, 21101 Novi Sad, Serbia, {\tt dockie\,@\,dmi.uns.ac.rs}.  The first named author gratefully acknowledges the support of Grant No.~174019 of the Ministry of Education, Science, and Technological Development of the Republic of Serbia, and Grant No.~0851/2015 of the Secretariat of Science and Technological Development of the Autonomous Province of Vojvodina.}
\ \  and \ 
James East%
\footnote{Centre for Research in Mathematics; School of Computing, Engineering and Mathematics, Western Sydney University, Locked Bag 1797, Penrith NSW 2751, Australia, {\tt J.East\,@\,WesternSydney.edu.au}.}%
}

\maketitle

\vspace{-0.5cm}

\begin{abstract}
We investigate the structure of the twisted Brauer monoid $\Bnt$, comparing and contrasting it to the structure of the (untwisted) Brauer monoid $\B_n$.  We characterise Green's relations and pre-orders on $\Bnt$, describe the lattice of ideals, and give necessary and sufficient conditions for an ideal to be idempotent-generated.  We obtain formulae for the rank (smallest size of a generating set) and (where applicable) the idempotent rank (smallest size of an idempotent generating set) of each principal ideal; in particular, when an ideal is idempotent-generated, its rank and idempotent rank are equal.  As an application of our results, we also describe the idempotent-generated subsemigroup of $\Bnt$ (which is not an ideal) as well as the singular ideal of $\Bnt$ (which is neither principal nor idempotent-generated), and we deduce a result of Maltcev and Mazorchuk that the singular part of the Brauer monoid $\B_n$ is idempotent-generated.

{\it Keywords}: Brauer monoids, twisted Brauer monoids, idempotents, ideals, rank, idempotent rank.

MSC: 20M20; 20M17, 05A18, 05E15.
\end{abstract}

\section{Introduction}\label{sect:intro}

The \emph{Temperley-Lieb algebras} were introduced in \cite{TL1971} to study lattice problems in (planar) statistical mechanics.  These algebras have played important roles in many different areas of mathematics, most notably in foundational works of Jones \cite{Jones1987} and Kauffman \cite{Kauffman1990} on knot polynomials.  As noted by Kauffman in \cite{Kauffman1990}, the structure of the Temperley-Lieb algebra is governed by an underlying (countably infinite) monoid that has now become known as the \emph{Kauffman monoid} \cite{BDP2002,LF2006}; an approach via a natural finite quotient of this monoid was described in \cite{Wilcox2007}.  
%
%
Kauffman also noted in \cite{Kauffman1990} that the Temperley-Lieb algebras are closely related to the algebras introduced by Brauer in his famous 1937 article \cite{Brauer1937} on invariant theory and representations of orthogonal groups.  
The Temperley-Lieb and Brauer algebras both have bases consisting of certain diagrams that are concatenated in a natural way (see below), so that the product of two basis elements is a scalar multiple of another basis element.
%
Other such algebras, known collectively as \emph{diagram algebras}, include partition algebras \cite{HR2005,Martin1994,Jones1994_2},
partial Brauer algebras \cite{HD2014,MarMaz2014}, Motzkin algebras \cite{BH2014}, rook monoid algebras \cite{HR2001,Solomon2002}, and many more.  These diagram algebras are all \emph{twisted semigroup algebras} \cite{Wilcox2007} of certain finite \emph{diagram semigroups} (such as the partition monoid, Brauer monoid, and Jones monoid), but they may also be viewed as (ordinary) semigroup algebras of the so-called \emph{twisted diagram semigroups} (the Kauffman monoid is a canonical example).

Studies of diagram semigroups have led to important results concerning the associated algebras, including cellularity \cite{Wilcox2007}, presentations \cite{JEgrpm,JEpnsn} and idempotent enumeration \cite{DEEFHHL1,DEEFHHL2}; see also \cite{EastGray} for an alternative approach to calculating dimensions of irreducible representations.  But it is also interesting to note that diagram semigroups have played a part in the development of semigroup theory itself, particularly in the context of regular $*$-semigroups \cite{EF,JEpme} and pseudovarieties of finite semigroups \cite{ADV2012_2,Auinger2014,Auinger2012}.  
%
Although the \emph{twisted} diagram semigroups are more closely related to diagram algebras, they have so far received less attention than their untwisted relatives, with existing studies \cite{BDP2002,DP2003,LF2006,ACHLV2015,BL2005,DE_Kauffman,DEEFHHL2} focusing mostly on the Kauffman monoid (which we have already discussed).  This article therefore aims to further the study of twisted diagram semigroups, and here we focus on the \emph{twisted Brauer monoid}.\footnote{The twisted Brauer monoid also played a role in the article \cite{ACHLV2015}, where it was called the \emph{wire monoid}.  We use the current terminology because of the above-mentioned links with twisted semigroup algebras.}  
In particular, we conduct a thorough investigation of the algebraic structure of the monoid, paying particular attention to Green's relations and pre-orders (which govern divisibility in the monoid and formalise several natural parameters associated to Brauer diagrams) and the lattice of ideals (which plays an important role in the cellular structure of the associated algebra \cite{GL1996}).
%
We also consider combinatorial problems such as determining which ideals are idempotent-generated, and calculating invariants such as the smallest size of (idempotent) generating sets.



The article is organised as follows.  In Section \ref{sect:Bn}, we recall the definition of the Brauer monoid $\B_n$, and record some known results we will need in what follows.  Section \ref{sect:Bnt}, which concerns the twisted Brauer monoid $\Bnt$, forms the bulk of the article, and consists of four subsections.  In Section \ref{sect:green}, we describe Green's relations and pre-orders on $\Bnt$, and we also characterise the regular elements of $\Bnt$.  Section \ref{sect:ideals} contains a classification of the ideals of $\Bnt$.  We calculate the smallest size of a generating set for each principal ideal of $\Bnt$ in Section \ref{sect:rank}, where we also give necessary and sufficient conditions for an ideal to be idempotent-generated; we also calculate the smallest size of an idempotent generating set for such an ideal.  Finally, in Section \ref{sect:applications}, we apply the results of previous sections to prove results about the singular part of $\Bnt$ and the idempotent-generated subsemigroups of $\Bnt$ and $\B_n$.


\section{The Brauer monoid}\label{sect:Bn}


Fix a non-negative integer $n$, and write $\bn=\{1,\ldots,n\}$ and $\bn'=\{1',\ldots,n'\}$.   Denote by $\B_n$ the set of all set partitions of $\bn\cup\bn'$ into blocks of size $2$.  For example, here is an element of $\B_6$:
\[
\al = \big\{
\{1,3\}, \{2,3'\}, \{4,1'\}, \{5,6\}, \{2',6'\}, \{4',5'\}
\big\}.
\]
There is a unique element of $\B_0$, namely the empty partition.  It is easy to see that 
\[
|\B_n|=(2n-1)!!=(2n-1)\cdot(2n-3)\cdots3\cdot1=\frac{(2n)!}{2^n\cdot n!} = \frac{n!}{2^n}\cdot{2n\choose n}.
\]
An element of $\B_n$ may be represented (uniquely) by a graph on vertex set $\bn\cup\bn'$; a single edge is included between vertices $u,v\in\bn\cup\bn'$ if and only if $\{u,v\}$ is a block of $\al$.  Such a graph is called a \emph{Brauer $n$-diagram} (or just a \emph{Brauer diagram} if $n$ is understood from context).  We typically identify $\al\in\B_n$ with its corresponding Brauer diagram.  When drawing a Brauer diagram, the vertices $1,\ldots,n$ are arranged in a horizontal line, with vertices $1',\ldots,n'$ in a parallel line below; unless otherwise specified, the vertices are assumed to be increasing from left to right.  For example, with $\al\in\B_6$ as above, we have:
\[
\al = \custpartn{1,2,3,4,5,6}{1,2,3,4,5,6}{\uarcx13{.5}\uarc56\darcx45{.3}\stline23\stline41\darcx26{.6}}.
\]
It will often be convenient to order the top and/or bottom vertices differently, but the ordering will always be made clear (see Figure \ref{fig:alphabetagamma}, for example).

The set $\B_n$ forms a monoid, known as the \emph{Brauer monoid} of degree $n$, under an operation we now describe.  Let $\al,\be\in\B_n$.  Write $\bn''=\{1'',\ldots,n''\}$.  Let $\alb$ be the graph obtained from $\al$ by changing the label of each lower vertex $i'$ to~$i''$.  Similarly, let $\beb$ be the graph obtained from $\be$ by changing the label of each upper vertex~$i$ to~$i''$.  Consider now the graph $\Ga(\al,\be)$ on the vertex set~$\bn\cup \bn'\cup \bn''$ obtained by joining $\alb$ and~$\beb$ together so that each lower vertex $i''$ of $\alb$ is identified with the corresponding upper vertex $i''$ of $\beb$.  Note that $\Ga(\al,\be)$, which we call the \emph{product graph of $\al,\be$}, may contain parallel edges.  We define $\al\be\in\B_n$ to be the Brauer diagram that has an edge $\{x,y\}$ if and only if $x,y\in\bn\cup\bn'$ are connected by a path in $\Ga(\al,\be)$.  An example calculation (with $n=10$) is given in Figure~\ref{fig:multinB10}.  

\begin{figure}[h]
\begin{center}
\begin{tikzpicture}[xscale=.4,yscale=0.4]
\begin{scope}[shift={(0,0)}]	
\uverts{1,...,10}
\lverts{1,...,10}
\stlines{3/3,4/6,6/7,7/8}
\uarc12
\uarc9{10}
\uarcx58{.6}
\darc12
\darc45
\darc9{10}
\draw(0.5,1)node[left]{$\al=$};
\end{scope}
\begin{scope}[shift={(0,-4)}]	
\uverts{1,...,10}
\lverts{1,...,10}
\stlines{3/2,9/9}
\uarc24
\uarc67
\uarc8{10}
\uarcx15{.7}
\darc13
\darc45
\darc78
\darcx6{10}{.7}
\draw(0.5,1)node[left]{$\be=$};
\draw(13,3)node{$\longrightarrow$};
\end{scope}
\begin{scope}[shift={(15,-1)}]	
\uverts{1,...,10}
\lverts{1,...,10}
\stlines{3/3,4/6,6/7,7/8}
\uarc12
\uarc9{10}
\uarcx58{.6}
\darc12
\darc45
\darc9{10}
\end{scope}
\begin{scope}[shift={(15,-3)}]	
\uverts{1,...,10}
\lverts{1,...,10}
\stlines{3/2,9/9}
\uarc24
\uarc67
\uarc8{10}
\uarcx15{.7}
\darc13
\darc45
\darc78
\darcx6{10}{.7}
\draw(13,2)node{$\longrightarrow$};
\end{scope}
\begin{scope}[shift={(30,-2)}]	
\uverts{1,...,10}
\lverts{1,...,10}
\stlines{3/2,7/9}
\uarc12
\uarc9{10}
\uarcx46{.4}
\uarcx58{.6}
\darc13
\darc45
\darc78
\darcx6{10}{.7}
\draw(10.5,1)node[right]{$=\al\be$};
\end{scope}
\end{tikzpicture}
\end{center}
\vspace{-5mm}
\caption{Two Brauer diagrams $\al,\be\in\B_{10}$ (left), their product $\al\be\in\B_{10}$ (right), and the graph $\Ga(\al,\be)$ (centre).}
\label{fig:multinB10}
\end{figure}

The identity element of $\B_n$ is the Brauer diagram $1=\custpartn{1,2,5}{1,2,5}{\stlines{1/1,2/2,5/5}\udotted25\ldotted25}$.  The set
\[
\S_n = \set{\al\in\B_n}{\dom(\al)=\codom(\al)=\bn}
\]
is the group of units of $\B_n$, and is (isomorphic to) the symmetric group on $\bn$.

Let $\al\in\B_n$.  A block of $\al$ is called a \emph{transversal} if it has non-empty intersection with both $\bn$ and $\bn'$, and an \emph{upper hook} (resp., \emph{lower hook}) if it is contained in $\bn$ (resp., $\bn'$).  The \emph{rank} of $\al$, denoted $\rank(\al)$, is equal to the number of transversals of $\al$.  For $x\in\bn\cup\bn'$, let $[x]_\al$ denote the block of $\al$ containing~$x$.  We define the \emph{domain} and \emph{codomain} of $\al$ to be the sets
\begin{align*}
\dom(\al) = \bigset{ x\in\bn } { [x]_\al\cap\bn'\not=\emptyset} &\AND
\codom(\al) = \bigset{ x\in\bn } { [x']_\al\cap\bn\not=\emptyset}.
\intertext{Note that $\rank(\al)=|\dom(\al)|=|\codom(\al)|$, and that $n-\rank(\al)$ is equal to the number of hooks of $\al$ (half of which are upper hooks, and half lower).  We also define the \emph{kernel} and \emph{cokernel} of $\al$ to be the equivalences}
\ker(\al) = \bigset{(x,y)\in\bn\times\bn}{[x]_\al=[y]_\al} &\AND
\coker(\al) = \bigset{(x,y)\in\bn\times\bn}{[x']_\al=[y']_\al}.
\end{align*}
To illustrate these ideas, with $\al = \custpartn{1,2,3,4,5,6}{1,2,3,4,5,6}{\uarcx13{.5}\uarc56\darcx45{.3}\stline23\stline41\darcx26{.6}}\in\B_6$ as above, we have $\rank(\al)=2$ and
\[
\dom(\al)=\{2,4\} \COMMA \codom(\al)=\{1,3\} \COMMA \ker(\al)=(\ 1,3\ | \ 2 \ |\ 4\ |\ 5,6\ ) \COMMA \coker(\al) = (\ 1\ | \ 2,6 \ |\ 3\ |\ 4,5\ ),
\]
using an obvious notation for equivalences.

It is immediate from the definitions that
\[
\begin{array}{rclcrcl}
\dom(\al\be) \hspace{-.25cm}&\sub&\hspace{-.25cm} \dom(\al), & &
\ker(\al\be)\hspace{-.25cm} &\sp&\hspace{-.25cm} \ker(\al),\\
\codom(\al\be) \hspace{-.25cm}&\sub&\hspace{-.25cm} \codom(\be), & &
\coker(\al\be)\hspace{-.25cm} &\sp&\hspace{-.25cm} \coker(\be),
\end{array}
\]
for all $\al,\be\in\B_n$.  For example, the identity $\ker(\al\be)\supseteq\ker(\al)$ says that any upper hook of $\al$ is an upper hook of $\al\be$. 

We now recall from \cite{EF} another way to specify an element of $\B_n$.  With this in mind, let $\al\in\B_n$.  We write
\begin{equation}
\label{eq:al}\tag{$\dagger$}
\al = \left( \begin{array}{c|c|c|c|c|c} 
i_1 \ & \ \cdots \ & \ i_r \ & \ a_1,b_1 \ & \ \cdots \ & \ a_s,b_s\ \ \\ \cline{4-6}
j_1 \ & \ \cdots \ & \ j_r \ & \ c_1,d_1 \ & \ \cdots \ & \ c_s,d_s \ \
\end{array} \!\!\! \right)
\end{equation}
to indicate that $\al$ has transversals $\{i_1,j_1'\},\ldots,\{i_r,j_r'\}$, upper hooks $\{a_1,b_1\},\ldots,\{a_s,b_s\}$, and lower hooks $\{c_1',d_1'\},\ldots,\{c_s',d_s'\}$.  Note that it is possible for either of $r,s$ to be $0$, but we always have $n=r+2s$.  In particular, we always have $\rank(\al)=r\equiv n\pmod2$.  
%

For $\al\in\B_n$, we write $\al^*$ for the Brauer diagram obtained from $\al$ by interchanging dashed and undashed vertices (i.e., by reflecting $\al$ in a horizontal axis).  
The ${}^*$ operation gives $\B_n$ the structure of a \emph{regular $*$-semigroup} \cite{NS1978}; that is, for all $\al,\be\in\B_n$,
\[
\al^{**}=\al \COMMA (\al\be)^*=\be^*\al^* \COMMA \al\al^*\al=\al \COMMA \al^*\al\al^*=\al^*.
\]
(The fourth identity follows quickly from the first three.)  This symmetry allows us to shorten many proofs.

Recall that \emph{Green's relations} $\gR,\gL,\gJ,\gH,\gD$ are defined on a semigroup $S$, for $x,y\in S$, by
\begin{gather*}
x\R y \iff xS^1=yS^1 \COMMA x\L y \iff S^1x=S^1y \COMMA x\J y \iff S^1xS^1 =S^1yS^1  , \\
{\H}={\R}\cap{\L} \COMMA {\D}={\R}\vee{\L}={\R}\circ{\L}={\L}\circ{\R}.
\end{gather*}
Here, $S^1$ denotes the monoid obtained by adjoining an identity $1$ to $S$ (if necessary).  If $S$ is finite, then ${\J}={\D}$.  
If $x\in S$, and if $\gK$ is one of Green's relations, we denote by $K_x$ the $\K$-class of $x$ in~$S$.  An $\H$-class contains an idempotent if and only if it is a group, in which case it is a maximal subgroup of $S$.  If $e$ and $f$ are $\D$-related idempotents of $S$, then the subgroups $H_e$ and $H_f$ are isomorphic.  If $S$ is a monoid, then the $\H$-class of the identity element of $S$ is the group of units of $S$.  
An element $x\in S$ is \emph{regular} if $x=xyx$ and $y=yxy$ for some $y\in S$ or, equivalently, if $D_x$ contains an idempotent, in which case $R_x$ and $L_x$ do, too.  In a $\D$-class of $S$, either every element is regular or every element is non-regular.  We say~$S$ is regular if every element of $S$ is regular.   For more background on semigroups, see for example \cite{Hig,Howie}.  The Brauer monoid $\B_n$ is regular since, as noted above, it is a regular $*$-semigroup.

The next result, which describes Green's relations on $\B_n$, was originally proved in \cite[Theorem 7]{Maz1998}; see also \cite{Larsson,Wilcox2007,FL2011}.

\ms
\begin{prop}[{Marorchuk \cite{Maz1998}}]\label{prop:Green_Bn}
Let $\al,\be\in\B_n$.  Then
\bit
\itemit{i} $\alpha \R \beta \iff \ker(\alpha)=\ker(\beta) \iff \al\S_n=\be\S_n$,
\itemit{ii} $\alpha \L \beta\iff\coker(\alpha)=\coker(\beta) \iff \S_n\al=\S_n\be$,
\itemit{iii} $\alpha \J \beta\iff\alpha \D \beta\iff\rank(\alpha) = \rank(\beta)\iff\S_n\al\S_n=\S_n\be\S_n$.  
\eit
In particular, 
\[\epfreseq
R_\al = \al\S_n \COMMA L_\al=\S_n\al \COMMA H_\al=\al\S_n\cap\S_n\al \COMMA D_\al=J_\al=\S_n\al\S_n  \qquad\text{for all $\al\in\B_n$.}
\]
\end{prop}

For the remainder of the paper, it will be convenient to define $z\in\{0,1\}$ with $z\equiv n\pmod2$.
We will also define the indexing set $I(n)=\{z,z+2,\ldots,n-2,n\}$.  So $\rank(\al)\in I(n)$ for all $\al\in\B_n$, and the $\D$-classes of~$\B_n$ are precisely the sets
\[
D_r = \set{\al\in\B_n}{\rank(\al)=r} \qquad\text{for $r\in I(n)$.}
\]
The following two results were proved in \cite[Theorem 8.4]{EastGray}.

\ms
\begin{prop}[East and Gray \cite{EastGray}]\label{prop:combinatorics_Bn}
Let $r=n-2s\in I(n)$, and put
\[
\rho_{nr}={n\choose r}\cdot(n-r-1)!!=\frac{n!}{2^ss!r!} 
\AND
\de_{nr}=\rho_{nr}^2\cdot r!=\frac{n!^2}{2^{2s}s!^2r!}.
\]
Then
\bit
\itemit{i} $D_r$ contains $\rho_{nr}$ $\R$-classes and $\rho_{nr}$ $\L$-classes, 
\itemit{ii} each $\H$-class contained in $D_r$ has size $r!$ (and group $\H$-classes contained in $D_r$ are isomorphic to~$\S_r$),
\itemit{iii} $|D_r|=\de_{nr}$.  \epf
\eit
\end{prop}


\ms
\begin{thm}[East and Gray \cite{EastGray}]\label{thm:ideals_Bn}
The ideals of $\B_n$ are precisely the sets
\[
I_r=D_z\cup D_{z+2}\cup\cdots\cup D_r = \set{\al\in\B_n}{\rank(\al)\leq r} \qquad\text{for $r\in I(n)$.}
\]
If $r\in I(n)\sm\{n\}$, then
\[
I_r=\la D_r\ra=\la E(D_r)\ra \AND \rank(I_r)=\idrank(I_r)=\rho_{nr},
\]
where the numbers $\rho_{nr}$ are defined in Proposition \ref{prop:combinatorics_Bn}. \epfres
\end{thm}

\section{The twisted Brauer monoid}\label{sect:Bnt}



When forming the product $\al\be$, where $\al,\be\in\B_n$, the product graph $\Ga(\al,\be)$ may contain components that lie completely in $\bn''$; such components are called \emph{floating components}.  We write $\tau(\al,\be)$ for the number of such floating components in $\Ga(\al,\be)$.  In the example from Figure \ref{fig:multinB10}, $\Ga(\al,\be)$ has a unique floating component, namely $\{1'',2'',4'',5''\}$, so $\tau(\al,\be)=1$.  
There are two main ways to modify the product in $\B_n$ to take these floating components into account.  One leads to the \emph{Brauer algebra} \cite{Brauer1937}, an associative algebra with $\B_n$ as its basis,
and the other leads to the \emph{twisted Brauer monoid}, which we now describe.  Specifically, we define
\[
\Bnt = \N\times\B_n = \set{(i,\al)}{i\in\N,\ \al\in\B_n}
\]
with product $\star$ defined, for $i,j\in\N$ and $\al,\be\in\B_n$, by
\[
(i,\al)\star(j,\be) = (i+j+\tau(\al,\be),\al\be).
\]
One easily checks that
\begin{equation}\label{eq:tau_identity}
\tau(\al,\be)+\tau(\al\be,\ga)=\tau(\al,\be\ga)+\tau(\be,\ga) \qquad\text{for all $\al,\be,\ga\in\B_n$.}
\end{equation}
It quickly follows that $\star$ is associative.  We call $\Bnt$ (with the $\star$ operation) the \emph{twisted Brauer monoid} of degree~$n$.  We note that there is a natural inclusion $\iota:\B_n\to\Bnt:\al\mt(0,\al)$, and we typically identify $\B_n$ with its image under $\iota$.  But it is important to note that $\iota$ is not a homomorphism, since $\al\star\be=(\tau(\al,\be),\al\be)\not=\al\be$ if $\tau(\al,\be)\not=0$.  It follows from the associativity of $\star$ that for any $\al_1,\ldots,\al_k\in\B_n$,
\[
\al_1\star\cdots\star\al_k=(\tau(\al_1,\ldots,\al_k),\al_1\cdots\al_k)
\]
for some $\tau(\al_1,\ldots,\al_k)\in\N$.  Note that for any $\al,\be,\ga\in\B_n$, $\tau(\al,\be,\ga)$ is equal to the common value in \eqref{eq:tau_identity}.  It is of special importance (and easily seen) that $\tau(\al,\be)=0$ if either of $\al,\be$ belongs to $\S_n$.  It is also immediate that $\tau(\al,\be)=\tau(\be^*,\al^*)$, so if we define $(i,\al)^*=(i,\al^*)$, then
\[
(i,\al)^{**}=(i,\al) \AND \big((i,\al)\star(j,\be)\big)^*=(j,\be)^*\star(i,\al)^*
\]
for all $(i,\al),(j,\be)\in\Bnt$.  In other words, $\Bnt$ is a $*$-semigroup (a semigroup with involution).  But this ${}^*$ operation does not give $\Bnt$ the structure of a regular $*$-semigroup \cite{NS1978}, since it is not necessarily the case that $(i,\al)\star(i,\al)^*\star(i,\al)=(i,\al)$; for example, this does not hold if $i\geq1$ or if $\tau(\al,\al^*)\geq1$.  In fact, $\Bnt$ is not a regular $*$-semigroup at all, as it is not even regular, as we will see in the next section.  

\subsection{Green's relations and pre-orders}\label{sect:green}

Our next goal is to describe Green's relations on the twisted Brauer monoid $\Bnt$.  In order to do this, it will be convenient to first describe Green's pre-orders on $\Bnt$.
Recall that Green's pre-orders $\leq_{\R},\leq_{\L},\leq_{\J}$ are defined on a semigroup $S$, for $x,y\in S$, by
\[
x\leqR y \iff xS^1\sub yS^1 \COMMA
x\leqL y \iff S^1x\sub S^1y \COMMA
x\leqJ y \iff S^1xS^1\sub S^1yS^1.
\]
So, for example, ${\R}={\leqR}\cap{\geqR}$.  In order to avoid confusion, we will use the symbols $\R$, $\leqR$, etc., for Green's relations and pre-orders on $\B_n$, and write $\Rt$, $\leqRt$, etc., for the corresponding relations and pre-orders on~$\Bnt$.  We first need to prove a result concerning Green's pre-orders on $\B_n$, which involves the twisting map $\tau$.


\ms
\begin{prop}\label{prop:Green_Bn_pre-order}
Let $\al,\be\in\B_n$.  Then
\bit
\itemit{i} $\alpha \leqR \beta \iff \ker(\alpha)\supseteq\ker(\beta)\iff \al=\be\de$ for some $\de\in\B_n$ with $\tau(\be,\de)=0$,
\itemit{ii} $\alpha \leqL \beta \iff \coker(\alpha)\supseteq\coker(\beta)\iff \al=\ga\be$ for some $\ga\in\B_n$ with $\tau(\ga,\be)=0$,
\itemit{iii} $\alpha \leqJ \beta \iff \rank(\alpha)\leq\rank(\beta)\iff \al=\ga\be\de$ for some $\ga,\de\in\B_n$ with $\tau(\ga,\be,\de)=0$.
\eit
\end{prop}

\pf We begin with (i).  Again, it is well-known that $\alpha \leqR \beta \iff \ker(\alpha)\supseteq\ker(\beta)$.  Next, suppose $\ker(\al)\supseteq\ker(\be)$.  Then we may write
\begin{align*}
\al &= \left( \begin{array}{c|c|c|c|c|c|c|c|c} 
i_1 \ & \ \cdots \ & \ i_r \ & \ a_1,b_1 \ & \ \cdots \ & \ a_s,b_s \ & \ a_{s+1},b_{s+1} \ & \ \cdots \ & \ a_{s+t},b_{s+t}  \ \ \\ \cline{4-9}
j_1 \ & \ \cdots \ & \ j_r \ & \ c_1,d_1 \ & \ \cdots \ & \ c_s,d_s \ & \ c_{s+1},d_{s+1} \ & \ \cdots \ & \ c_{s+t},d_{s+t}  \ \  
\end{array} \!\!\! \right)
\intertext{and}
\be &= \left( \begin{array}{c|c|c|c|c|c|c|c|c|c|c} 
i_1 \ & \ \cdots \ & \ i_r \ & \ a_1 \ & \ b_1 \ & \ \cdots \ & \ a_s \ & \ b_s \ & \ a_{s+1},b_{s+1} \ & \ \cdots \ & \ a_{s+t},b_{s+t}  \ \ \\ \cline{9-11}
k_1 \ & \ \cdots \ & \ k_r \ & \ e_1 \ & \ f_1 \ & \ \cdots \ & \ e_s \ & \ f_s \ & \ e_{s+1},f_{s+1} \ & \ \cdots \ & \ e_{s+t},f_{s+t}  \ \  
\end{array} \!\!\! \right).
\intertext{It is easy to check that $\al=\be\de$ with $\tau(\be,\de)=0$, where}
\de &= \left( \begin{array}{c|c|c|c|c|c|c|c|c|c|c} 
k_1 \ & \ \cdots \ & \ k_r \ & \ e_{s+1} \ & \ f_{s+1} \ & \ \cdots \ & \ e_{s+t} \ & \ f_{s+t} \ & \ e_1,f_1 \ & \ \cdots \ & \ e_s,f_s  \ \ \\ \cline{9-11}
j_1 \ & \ \cdots \ & \ j_r \ & \ c_{s+1} \ & \ d_{s+1} \ & \ \cdots \ & \ c_{s+t} \ & \ d_{s+t} \ & \ c_1,d_1 \ & \ \cdots \ & \ c_s,d_s  \ \  
\end{array} \!\!\! \right).
\end{align*}
This completes the proof of (i).  Part (ii) follows by duality.  

For (iii), suppose $\rank(\al)\leq\rank(\be)$.  As above, it suffices to demonstrate the existence of $\ga,\de$ with the desired properties.  We may write
\begin{align*}
\al &= \left( \begin{array}{c|c|c|c|c|c|c|c|c} 
i_1 \ & \ \cdots \ & \ i_r \ & \ a_1,b_1 \ & \ \cdots \ & \ a_s,b_s \ & \ a_{s+1},b_{s+1} \ & \ \cdots \ & \ a_{s+t},b_{s+t}  \ \ \\ \cline{4-9}
j_1 \ & \ \cdots \ & \ j_r \ & \ c_1,d_1 \ & \ \cdots \ & \ c_s,d_s \ & \ c_{s+1},d_{s+1} \ & \ \cdots \ & \ c_{s+t},d_{s+t}  \ \  
\end{array} \!\!\! \right)
\intertext{and}
\be &= \left( \begin{array}{c|c|c|c|c|c|c|c|c|c|c} 
l_1 \ & \ \cdots \ & \ l_r \ & \ g_1 \ & \ h_1 \ & \ \cdots \ & \ g_s \ & \ h_s \ & \ g_{s+1},h_{s+1} \ & \ \cdots \ & \ g_{s+t},h_{s+t}  \ \ \\ \cline{9-11}
k_1 \ & \ \cdots \ & \ k_r \ & \ e_1 \ & \ f_1 \ & \ \cdots \ & \ e_s \ & \ f_s \ & \ e_{s+1},f_{s+1} \ & \ \cdots \ & \ e_{s+t},f_{s+t}  \ \  
\end{array} \!\!\! \right).
\intertext{Now put}
\ve &= \left( \begin{array}{c|c|c|c|c|c|c|c|c} 
l_1 \ & \ \cdots \ & \ l_r \ & \ g_1,h_1 \ & \ \cdots \ & \ g_s,h_s \ & \ g_{s+1},h_{s+1} \ & \ \cdots \ & \ g_{s+t},h_{s+t}  \ \ \\ \cline{4-9}
j_1 \ & \ \cdots \ & \ j_r \ & \ c_1,d_1 \ & \ \cdots \ & \ c_s,d_s \ & \ c_{s+1},d_{s+1} \ & \ \cdots \ & \ c_{s+t},d_{s+t}  \ \  
\end{array} \!\!\! \right).
\end{align*}
Then $\ker(\ve)\supseteq\ker(\be)$.  By (i), it follows that there exists $\de\in\B_n$ with $\ve=\be\de$ and $\tau(\be,\de)=0$.  But also $\ve\L\al$, so Proposition \ref{prop:Green_Bn}(ii) gives $\al=\ga\ve$ for some $\ga\in\S_n$.  In particular, $\al=\ga\ve=\ga\be\de$, and $\tau(\ga,\be,\de)=\tau(\ga,\be\de)+\tau(\be,\de)=0$. \epf


\ms
\begin{prop}\label{prop:Green_Bnt_pre-order}
Let $i,j\in\N$ and $\al,\be\in\B_n$.  If $\K$ is any of $\R,\L,\J$, then
\[
\text{$(i,\al)\leqKt(j,\be) \iff i\geq j$ and $\al\leqK\be$.}
\]
\end{prop}

\pf We just treat the $\leqJt$ pre-order, since the other cases are similar.  Suppose first that $(i,\al)\leqJt(j,\be)$.  Then
\[
(i,\al) = (h,\ga)\star(j,\be)\star(k,\de) = (h+j+k+\tau(\ga,\be,\de),\ga\be\de)
\]
for some $h,k\in\N$ and $\ga,\de\in\B_n$.  But then $i=h+j+k+\tau(\ga,\be,\de)\geq j$ and $\al=\ga\be\de\leqJ\be$.  Conversely, suppose $i\geq j$ and $\al\leqJ\be$.  By Proposition \ref{prop:Green_Bn_pre-order}(iii), there exists $\ga,\de\in\B_n$ such that $\al=\ga\be\de$ and $\tau(\ga,\be,\de)=0$.  But then one easily checks that $(i,\al)= (i-j,\ga)\star(j,\be)\star(0,\de)$, completing the proof. \epf


Let $i\in\N$ and $\al\in\B_n$.  If $\K$ is one of Green's relations, we write $K_\al$ and $K_{(i,\al)}^\tau$ for the $\K$-class of $\al$ in $\B_n$ and the $\Kt$-class of $(i,\al)$ in $\Bnt$.

\ms
\begin{cor}\label{cor:Green_Bnt}
Let $i,j\in\N$ and $\al,\be\in\B_n$.  If $\K$ is any of ${\R},{\L},{\H},{\J},{\D}$, then
\[
\text{\emph{$(i,\al)\Kt(j,\be)\iff i=j$ and $\al\K\be$.}}
\]
Consequently, $K_{(i,\al)}^\tau=\{i\}\times K_\al$ for any $(i,\al)\in\Bnt$. 
\end{cor}

\pf The descriptions of the $\Rt,\Lt,\Ht,\Jt$ relations follow immediately from Proposition \ref{prop:Green_Bnt_pre-order}.  It remains only to show that 
${\Jt}\sub{\Dt}$.  But this is true because
\begin{align*}
(i,\al)\J(j,\be) \implies [\text{$i=j$ and $\al\J\be$} ]
&\implies [\text{$i=j$ and $\al\D\be$}] 
\implies [\text{$i=j$ and $\al\R\ga\L\be$} \ \ \text{for some $\ga\in\B_n$} ]\\
\epfreseq
&\implies (i,\al)\Rt(i,\ga)\Lt(j,\be) 
\implies (i,\al)\Dt(j,\be).
\end{align*}




%

So the $\Dt$-classes of $\Bnt$ are precisely the sets
\[
\DClass rk = \{k\}\times D_r = \set{(k,\al)}{\rank(\al)=r} \qquad\text{for $r\in I(n)$ and $k\in\N$.}
\]
Note that under the identification of $\al\in\B_n$ with $(0,\al)\in\Bnt$, we have $\DClass r0=D_r$ for all $r\in I(n)$.

Recall that the set $S/{\J}$ of all $\J$-classes of a semigroup $S$ is a partially ordered set under the order $\leq$ defined, for $x,y\in S$, by $J_x\leq J_y \iff x\leqJ y$.  We will write $\leq$ and $\leqt$ for the partial orders on $\B_n/{\D}$ and $\Bnt/{\Dt}$, respectively (recall that ${\J}={\D}$ and ${\Jt}={\Dt}$ in $\B_n$ and $\Bnt$).  So, by Propositions \ref{prop:Green_Bn_pre-order} and \ref{prop:Green_Bnt_pre-order}, we have
\[
D_r\leq D_s \iff r\leq s \AND \DClass rk \leq \DClass sl \iff [\text{$r\leq s$ and $k\geq l$}].
\]
So the partially ordered set $(\Bnt/{\Dt},{\leqt})$ is a lattice, and is order-isomorphic to the direct product of the chains $(I(n),\leq)$ and $(\N,\geq)$; this is analogous to the case of the Kauffman monoid \cite{LF2006}.  
Figure \ref{fig:partial_order} gives an illustration for $n=7$ (the reader may ignore the shading in the diagram for now).  

\begin{figure}[h]
\begin{center}
\begin{tikzpicture}[scale=0.9]
\draw[gray!40,rounded corners, fill=gray!40] (15+1,-5-1-.33)--(6-1,-2-1+.33)--(6-1,2+1+.33)--(15+1,-1+1-.33);
\draw [gray!60,rounded corners, fill=gray!60] (6-.75,0)--(6-.75,-2-.75+.25)--(9+.75,-3-.75-.25)--(9+.75,1+.75-.25)--(6-.75,2+.75+.25)--(6-.75,0); 
\DRow6
\DRow4
\DRow2
\DRow0
\DColumn{0}{0}{0}
\DColumn{3}{-1}{1}
\DColumn{6}{-2}{2}
\DColumn{9}{-3}{3}
\DColumn{12}{-4}{4}
\end{tikzpicture}
\end{center}
\vspace{-5mm}
\caption{The structure of the partially ordered set $(\B_7^\tau/{\Dt},{\leqt})$.  The principal ideal $\Ideal52$ is shaded light grey, and its generating set $\Minimal52$ is shaded dark grey.}
\label{fig:partial_order}
\end{figure} 


We conclude this section with a description of the regular elements of $\Bnt$.

\ms
\begin{prop}\label{prop:regular}
An element $(i,\al)\in\Bnt$ is regular if and only if $i=0$ and $\rank(\al)>0$.  In particular, $\Bnt$ is not regular.
\end{prop}

\pf 
From $(i,\al)\star(j,\be)\star(i,\al)=(2i+j+\tau(\al,\be,\al),\al\be\al)$, we deduce that $(i,\al)$ cannot be regular  if
\bit
\item[(i)]  $i\geq1$, since then $2i+j+\tau(\al,\be,\al)\geq2i>i$, or if
\item[(ii)]  $\rank(\al)=0$, since then $2i+j+\tau(\al,\be,\al)\geq2i+j+1>i$. 
\eit
Conversely, if $i=0$ and $r>0$, then one may easily check that
$
\ve 
= \custpartn{1,3,4,5,6,8,9}{1,3,4,5,7,8,9}{\stlines{1/1,3/3}\uarcs{5/6,8/9}\udotteds{1/3,6/8}\ldotteds{1/3,5/7}\darcs{4/5,7/8}\uarcxhalf461\darcxhalfr691}\in\DClass r0 
$
is an idempotent of $\Bnt$ (i.e., $\ve=\ve\star\ve$).  It follows that the $\Dt$-classes $\DClass r0$ with $r>0$ are all regular. \epf

\subsection{Ideals}\label{sect:ideals}

We may now describe the ideals of $\Bnt$.  Recall that a \emph{principal ideal} of a semigroup $S$ is of the form
\[
S^1aS^1=\set{xay}{x,y\in S^1} = \set{x\in S}{x\leqJ a} \qquad\text{for $a\in S$.}
\]
By Proposition \ref{prop:Green_Bnt_pre-order}, we may immediately describe the principal ideals of $\Bnt$.  These are precisely the sets
\[
\Ideal rk 
= \set{(i,\al)}{\rank(\al)\leq r,\ i\geq k} \qquad\text{for $r\in I(n)$ and $k\in\N$.}
\]
Note that $\Ideal rk\sub\Ideal sl \iff \DClass rk\leqt\DClass sl \iff [\text{$r\leq s$ and $k\geq l$}]$.
The principal ideal $\Ideal 52$ of $\B_7^\tau$ is pictured in Figure \ref{fig:partial_order}.  
We now show that every ideal of $\Bnt$ is the union of finitely many principal ideals (not every infinite semigroup shares this property).

\ms\ms\ms\ms
\begin{prop}
\bit
\itemit{i} Let $r_1,\ldots,r_s\in I(n)$ and $k_1,\ldots,k_s\in\N$, with $r_1>\cdots>r_s$ and $k_1>\cdots>k_s$.  Then
$
\Ideal{r_1}{k_1} \cup\cdots\cup\Ideal{r_s}{k_s}
$
is an ideal of $\Bnt$.  
\itemit{ii} Each ideal of $\Bnt$ is of the form described in \emph{(i)}.
\itemit{iii} Each ideal of $\Bnt$ is uniquely determined by (and uniquely determines) the parameters $r_1,\ldots,r_s,k_1,\ldots,k_s$, as described in \emph{(i)}.
\eit
\end{prop}

\pf Part (i) is clear.
Next, suppose $I$ is an arbitrary non-empty ideal of $\Bnt$.  Put
\[
r_1=\max\set{\rank(\al)}{(k,\al)\in I\ (\exists k\in\N)} \AND k_1=\min\set{k\in\N}{(k,\al)\in I \ (\exists \al\in D_{r_1})}.
\]
(Note that $k_1$ is defined in terms of $r_1$.)  Then $\Ideal{r_1}{k_1}\sub I$.  If $I=\Ideal{r_1}{k_1}$, then we are done.  Otherwise, put
\[
r_2=\max\set{\rank(\al)}{(k,\al)\in I\sm\Ideal{r_1}{k_1}\ (\exists k\in\N)} \quad\text{and}\quad k_2=\min\set{k\in\N}{(k,\al)\in I\sm\Ideal{r_1}{k_1} \ (\exists \al\in D_{r_2})}.
\]
(Note that $r_1>r_2$ is obvious, while $k_1>k_2$ follows from the fact that $I_{r_1;k_1}$ already contains $I_{r_2;k_1}$.)
Then $\Ideal{r_2}{k_2}\sub I$.  If $I=\Ideal{r_1}{k_1}\cup\Ideal{r_2}{k_2}$, then we are done.  Otherwise, we similarly define $r_3$ and $k_3$.  Continuing in this fashion, 
since $I(n)$ is a finite chain, we eventually obtain
$
I = \Ideal{r_1}{k_1} \cup\cdots\cup\Ideal{r_s}{k_s}
$
for some $r_1,\ldots,r_s\in I(n)$ and $k_1,\ldots,k_s\in\N$ with $r_1>\cdots>r_s$ and $k_1>\cdots>k_s$, giving (ii).
For (iii), it is clear that $\Ideal{r_1}{k_1} \cup\cdots\cup\Ideal{r_s}{k_s}=\Ideal{q_1}{l_1} \cup\cdots\cup\Ideal{q_t}{l_t}$ if and only if $(r_1,\ldots,r_s)=(q_1,\ldots,q_t)$ and $(k_1,\ldots,k_s)=(l_1,\ldots,l_t)$.~\epf

\ms
\begin{rem}
%
Note that
\begin{align*}
\Ideal{r_1}{k_1} \cup\cdots\cup\Ideal{r_s}{k_s}\sub\Ideal{q_1}{l_1} \cup\cdots\cup\Ideal{q_t}{l_t}
&\iff
(\forall i\in\bs)(\exists j\in\bt)\  \Ideal{r_i}{k_i}\sub\Ideal{q_j}{l_j}\\
&\iff
(\forall i\in\bs)(\exists j\in\bt)\  [\text{$r_i\leq q_j$ and $k_i\geq l_j$}].
\end{align*}
\end{rem}



\subsection{Small generating sets}\label{sect:rank}

We now turn to the question of minimal generation of the principal ideals.  Recall that if $S$ is a semigroup, then the \emph{rank} of $S$, denoted $\rank(S)$, is the minimum cardinality of a subset $A\sub S$ such that $S=\la A\ra$.  If $S$ is idempotent-generated, then the \emph{idempotent rank} of $S$, denoted $\idrank(S)$, is defined analogously with respect to generating sets consisting of idempotents.  In this section, we give necessary and sufficient conditions for a principal ideal $\Ideal rk$ to be idempotent generated.  We also calculate the rank and idempotent rank (if appropriate) for an arbitrary principal ideal $\Ideal rk$; in particular, we show that $\rank(\Ideal rk)=\idrank(\Ideal rk)$ if $\Ideal rk$ is idempotent-generated.

If $\Si\sub\B_n$ (resp., $\Ga\sub\Bnt$), we write $\la\Si\ra$ (resp., $\dla\Ga\dra$) for the subsemigroup of $\B_n$ (resp., $\Bnt$) generated by $\Si$ (resp., $\Ga$).  Since we identify $\B_n$ with a subset of $\Bnt$, via the mapping $\al\mt(0,\al)$, it is possible to consider both $\la\Si\ra$ and $\dla\Si\dra$ for a subset $\Si\sub\B_n$; these are obviously not equal in general.

It will be necessary to consider the ideals $\Ideal rk$ in a number of separate cases, depending on the values of the parameters $r,k$ (see Theorem \ref{thm:rank}).  We begin with the ideals $\Ideal r0$ with $r<n$.  For this, we will need the following two lemmas, the second of which will also be used later.

\ms
\begin{lemma}\label{lem:Dr_Dr+2}
If $r\leq n-4$, then $D_r\sub D_{r+2}\star D_{r+2}$.
\end{lemma}

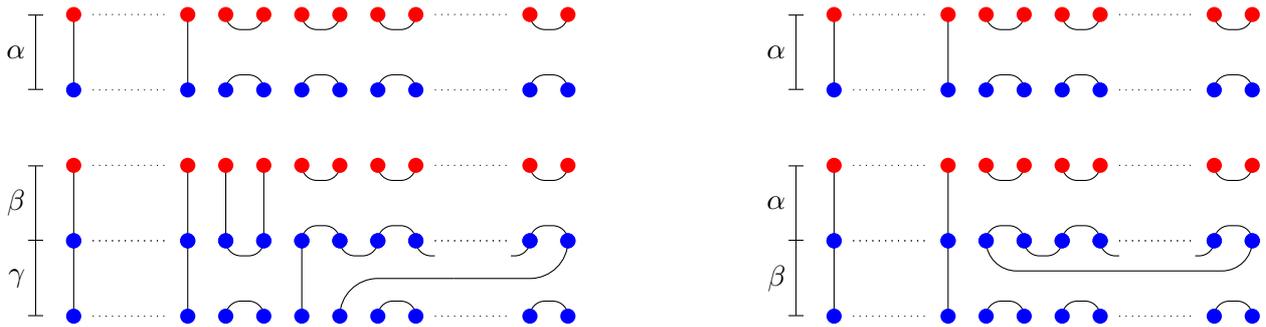
\begin{figure}[h]
\begin{center}
\begin{tikzpicture}[scale=.5]
\begin{scope}[shift={(0,0)}]	
\stlines{1/1,4/4}
\uarcs{5/6,7/8,9/10,13/14}
\darcs{5/6,7/8,9/10,13/14}
\udotteds{1/4,10/13}
\ldotteds{1/4,10/13}
\uvertcols{1,4,5,6,7,8,9,10,13,14}{red}
\lvertcols{1,4,5,6,7,8,9,10,13,14}{blue}
\draw[|-|] (0,0)--(0,2);
\draw(0,1)node[left]{$\al$};
\end{scope}
\begin{scope}[shift={(0,-4)}]	
\stlines{1/1,4/4,5/5,6/6}
\uarcs{7/8,9/10,13/14}
\darcs{7/8,9/10,13/14}
\udotteds{1/4,10/13}
\ldotteds{1/4,10/13}
\uvertcols{1,4,5,6,7,8,9,10,13,14}{red}
\lvertcols{1,4,5,6,7,8,9,10,13,14}{blue}
\draw[|-|] (0,0)--(0,2);
\draw(0,1)node[left]{$\be$};
\end{scope}
\begin{scope}[shift={(0,-6)}]	
\stlines{1/1,4/4,7/7}
\darcxhalf8{11}1
\uarcxhalfr{11}{14}1
\uarcs{5/6,8/9}
\uarchalf{10}{10.5}
\uarchalfr{12.5}{13}
\darcs{5/6,9/10,13/14}
\udotteds{1/4,10/13}
\ldotteds{1/4,10/13}
\uvertcols{1,4,5,6,7,8,9,10,13,14}{blue}
\lvertcols{1,4,5,6,7,8,9,10,13,14}{blue}
\draw[|-] (0,0)--(0,2);
\draw(0,1)node[left]{$\ga$};
\end{scope}
\begin{scope}[shift={(20,0)}]	
\stlines{1/1,4/4}
\uarcs{5/6,7/8,11/12}
\darcs{5/6,7/8,11/12}
\udotteds{1/4,8/11}
\ldotteds{1/4,8/11}
\uvertcols{1,4,5,6,7,8,11,12}{red}
\lvertcols{1,4,5,6,7,8,11,12}{blue}
\draw[|-|] (0,0)--(0,2);
\draw(0,1)node[left]{$\al$};
\end{scope}
\begin{scope}[shift={(20,-4)}]	
\stlines{1/1,4/4}
\uarcs{5/6,7/8,11/12}
\darcs{5/6,7/8,11/12}
\udotteds{1/4,8/11}
\ldotteds{1/4,8/11}
\uvertcols{1,4,5,6,7,8,11,12}{red}
\lvertcols{1,4,5,6,7,8,11,12}{blue}
\draw[|-|] (0,0)--(0,2);
\draw(0,1)node[left]{$\al$};
\end{scope}
\begin{scope}[shift={(20,-6)}]	
\stlines{1/1,4/4}
\uarcs{6/7}
\uarchalf8{8.5}
\uarchalfr{10.5}{11}
\uarcx5{12}{.8}
\darcs{5/6,7/8,11/12}
\udotteds{1/4,8/11}
\ldotteds{1/4,8/11}
\uvertcols{1,4,5,6,7,8,11,12}{blue}
\lvertcols{1,4,5,6,7,8,11,12}{blue}
\draw[|-] (0,0)--(0,2);
\draw(0,1)node[left]{$\be$};
\end{scope}
\end{tikzpicture}
\end{center}
\vspace{-5mm}
\caption{Diagrammatic verification that $\al=\be\ga$ with $\tau(\be,\ga)=0$ from the proof of Lemma \ref{lem:Dr_Dr+2} (left), and $\al=\al\be$ with $\tau(\al,\be)=1$ from the proof of Lemma \ref{lem:alphabeta} (right); see the text for more details.  In both cases, red vertices are ordered $i_1,\ldots,i_r,a_1,b_1,\ldots,a_s,b_s$, and blue vertices are ordered $j_1,\ldots,j_r,c_1,d_1,\ldots,c_s,d_s$.}
\label{fig:alphabetagamma}
\end{figure}

\pf Write $\al$ as in Equation \eqref{eq:al} on Page \pageref{eq:al},
where $r\leq n-4$.  We show in Figure \ref{fig:alphabetagamma} (left) that $\al=\be\ga$ for some $\be,\ga\in D_{r+2}$ with $\tau(\be,\ga)=0$.  \epf
%

\ms
\begin{rem}
A weaker version of Lemma \ref{lem:Dr_Dr+2} was proved in \cite[Lemma 8.3]{EastGray}, where it was shown that $D_r\sub\la D_{r+2}\ra$; the proof of that result was much simpler, as no conditions were imposed on the twisting map~$\tau$, and the $*$-regular structure of $\B_n$ played a role.
\end{rem}

\ms
\begin{lemma}\label{lem:alphabeta}
If $\al\in\BnSn$, then $\al=\al\be$ for some $\be\in\B_n$ with $\rank(\be)=\rank(\al)$ and $\tau(\al,\be)=1$.
\end{lemma}


\pf Write $\al$ as in Equation \eqref{eq:al} on Page \pageref{eq:al}.
We demonstrate the existence of $\be$ in Figure \ref{fig:alphabetagamma} (right). \epf


\ms
\begin{prop}\label{prop:Ir0=<Dr>}
If $r\in I(n)\sm\{n\}$, then $\Ideal r0=\dla D_r\dra$.
\end{prop}

\pf We first show, by descending induction, that $D_s\sub\dla D_r\dra$ for all $s\in I(r)$.  Indeed, this is obvious if $s=r$, while if $s<r$, then Lemma \ref{lem:Dr_Dr+2} and an induction hypothesis gives $D_s\sub D_{s+2}\star D_{s+2}\sub\dla D_r\dra$.  It follows that $I_r\sub\dla D_r\dra$.  Now suppose $i\geq1$ and $\al\in I_r$.  We have seen that $\al\in\dla D_r\dra$.  By Lemma \ref{lem:alphabeta}, we may choose some $\be\in D_r$ such that $\al=\al\be$ and $\tau(\al,\be)=1$.  But then it quickly follows that
\[
(i,\al)=\al\star\underbrace{\be\star\cdots\star\be}_{i}\in\dla D_r\dra.
\]
We have shown that $\Ideal r0\sub\dla D_r\dra$.  The reverse inclusion is clear. \epf

Proposition \ref{prop:Ir0=<Dr>} does not hold for the top ideal $\Ideal n0=\Bnt$, but we may use it as a stepping stone to calculate $\rank(\Bnt)$.  Recall that $\rank(\S_n)=2$ if $n\geq3$.  

\ms
\begin{prop}
Suppose $n\geq3$.  Let $\al,\be\in\S_n$ be such that $\S_n=\la\al,\be\ra$, and let $\ga\in \DClass{n-2}0$ and $(1,\de)\in\DClass n1$ be arbitrary.  Then $\Bnt=\dla\al,\be,\ga,(1,\de)\dra$.  Further, $\rank(\Bnt)=4$.
\end{prop}

\pf Write $S=\dla\al,\be,\ga,(1,\de)\dra$.  First note that $\S_n=\la\al,\be\ra=\dla\al,\be\dra\sub S$.  Together with Proposition~\ref{prop:Green_Bn}, it then follows that $S$ contains $D_{n-2}=D_\ga=\S_n\ga\S_n=\S_n\star\ga\star\S_n$.  Proposition \ref{prop:Ir0=<Dr>} then gives $\Ideal{n-2}0=\dla D_{n-2}\dra\sub S$.  Finally, let $i\geq1$ and $\si\in\S_n$ be arbitrary.  Then
\[
(i,\si) = (0,\si\de^{-i}) \star \underbrace{(1,\de)\star\cdots\star(1,\de)}_i\in S,
\]
completing the proof that $\Bnt=S=\dla\al,\be,\ga,(1,\de)\dra$.  It also follows that $\rank(\Bnt)\leq4$.  

Suppose now that $\Bnt=\dla\Si\dra$.  The proof will be complete if we can show that $|\Si|\geq4$.  Since $\BntSn=\Ideal{n-2}0\cup\Ideal n1$ is an ideal of $\Bnt$, it follows that $\Si$ contains a generating set for $\S_n$, so that $|\Si\cap\S_n|\geq2$.  Now let $\si\in\S_n$ be arbitrary, and consider an expression
\[
(1,\si)=(i_1,\al_1)\star\cdots\star(i_k,\al_k) = (i_1+\cdots+i_k+\tau(\al_1,\ldots,\al_k),\al_1\cdots\al_k),
\]
where $(i_1,\al_1),\ldots,(i_k,\al_j)\in\Si$.  Since $\al_1\cdots\al_k=\si\in\S_n$, and since $\BnSn$ is an ideal of $\B_n$, it follows that $\al_1,\ldots,\al_k\in\S_n$.  Then $\tau(\al_1,\ldots,\al_k)=0$, so $1=i_1+\cdots+i_k$, which gives $i_s=1$ for some (unique) $s\in\bk$.  Thus, $\Si$ contains an element of $\DClass n1$: namely, $(1,\al_s)$.  Similarly, consideration of an element of $\DClass {n-2}0$ as a product of elements from $\Si$ shows that $\Si$ contains an element of $\DClass{n-2}0$.  As noted above, this completes the proof.~\epf

Next, we calculate $\rank(\Ideal r0)$ in the case that $0<r<n$.  In fact, since the ideal $\Ideal r0$ is idempotent-generated for such a value of $r$ (as we will soon show), we will also calculate $\idrank(\Ideal r0)$.  
Since
\[
(i,\al)\star(i,\al)=(2i+\tau(\al,\al),\al^2),
\]
it follows that all idempotents of $\Bnt$ are contained in $\B_n$.  However, not every idempotent of $\B_n$ is an idempotent of $\Bnt$; that is, $\al=\al^2$ in $\B_n$ does not necessarily imply $\al=\al\star\al$ in $\Bnt$.    In order to avoid confusion when discussing idempotents from $\B_n$ and $\Bnt$, if $\Si\sub\B_n$, we will write 
\[
E(\Si) = \set{\al\in\Si}{\al=\al^2} \AND \Et(\Si) = \set{\al\in\Si}{\al=\al\star\al}.
\]
For example, one may easily check that
\[
\al=\custpartn{1,...,6}{1,...,6}{\stlines{1/1,3/3}\uarcs{2/4,5/6}\darcs{2/4,5/6}}\in E(\B_6)\sm\Et(\B_6^\tau) \qquad\text{but}\qquad
\be=\custpartn{1,...,6}{1,...,6}{\stlines{1/1,3/2}\uarcs{2/4,5/6}\darcx45{.3}\darcx36{.6}}\in \Et(\B_6^\tau)
\]
Indeed, $\al\star\al=(2,\al)\not=\al$ in $\B_6^\tau$.  The idempotents of $\B_n$ and $\Bnt$ (and a number of other diagram semigroups) were characterised and enumerated in \cite{DEEFHHL1}, but we will not need to use these descriptions here.

\ms
\begin{prop}\label{prop:Ir0=<E(Dr)>}
Suppose $r\in I(n)\sm\{0,n\}$.  Then $\Ideal r0=\dla \Et(D_r)\dra$.
\end{prop}

\begin{figure}[h]
\begin{center}
\begin{tikzpicture}[scale=.5]
\begin{scope}[shift={(0,0)}]	
\stlines{1/1,4/4,5/5,6/6}
\uarcs{7/8,9/10,13/14}
\darcs{7/8,9/10,13/14}
\udotteds{1/4,10/13}
\ldotteds{1/4,10/13}
\uvertcols{1,4,5,6,7,8,9,10,13,14}{red}
\lvertcols{1,4,5,6,7,8,9,10,13,14}{blue}
\draw[|-|] (0,0)--(0,2);
\draw(0,1)node[left]{$\al$};
\end{scope}
\begin{scope}[shift={(0,-2)}]	
\stlines{1/1,4/4,5/6,6/5,7/7,8/8,9/9,10/10,13/13,14/14}
\udotteds{1/4,10/13}
\ldotteds{1/4,10/13}
\uvertcols{1,4,5,6,7,8,9,10,13,14}{blue}
\lvertcols{1,4,5,6,7,8,9,10,13,14}{blue}
\draw[|-] (0,0)--(0,2);
\draw(0,1)node[left]{$\si_{ij}$};
\end{scope}
\begin{scope}[shift={(0,-6)}]	
\stlines{1/1,4/4,5/5,6/6}
\uarcs{7/8,9/10,13/14}
\darcs{7/8,9/10,13/14}
\udotteds{1/4,10/13}
\ldotteds{1/4,10/13}
\uvertcols{1,4,5,6,7,8,9,10,13,14}{red}
\lvertcols{1,4,5,6,7,8,9,10,13,14}{blue}
\draw[|-|] (0,0)--(0,2);
\draw(0,1)node[left]{$\al$};
\end{scope}
\begin{scope}[shift={(0,-8)}]	
\stlines{1/1,4/4,5/6,14/14}
\uarcs{6/7,8/9}
\uarchalf{10}{10.5}
\uarchalfr{12.5}{13}
\darchalf9{9.5}
\darchalfr{11.5}{12}
\darcs{7/8}
\darcx5{13}{.7}
\udotteds{1/4,10/13}
\ldotteds{1/4,9/12}
\uvertcols{1,4,5,6,7,8,9,10,13,14}{blue}
\lvertcols{1,4,5,6,7,8,9,12,13,14}{blue}
\draw[|-] (0,0)--(0,2);
\draw(0,1)node[left]{$\be_1$};
\end{scope}
\begin{scope}[shift={(0,-10)}]	
\stlines{1/1,4/4,6/6,7/5}
\uarcs{8/9,12/13}
\darcs{7/8,9/10,13/14}
\uarcx5{14}{.7}
\ldotteds{1/4,10/13}
\lvertcols{1,4,5,6,7,8,9,10,13,14}{blue}
\draw[|-] (0,0)--(0,2);
\draw(0,1)node[left]{$\be_2$};
\end{scope}
\begin{scope}[shift={(0,-15)}]	
\stlines{1/1,4/4,5/5}
\uarcs{6/7,8/9,12/13}
\darcs{6/7,8/9,12/13}
\udotteds{1/4,9/12}
\ldotteds{1/4,9/12}
\uvertcols{1,4,5,6,7,8,9,12,13}{red}
\lvertcols{1,4,5,6,7,8,9,12,13}{blue}
\draw[|-|] (0,0)--(0,2);
\draw(0,1)node[left]{$\al$};
\end{scope}
\begin{scope}[shift={(0,-17)}]	
\stlines{1/1,4/4,5/6,6/5,7/7,8/8,9/9,12/12,13/13}
\udotteds{1/4,9/12}
\ldotteds{1/4,9/12}
\uvertcols{1,4,5,6,7,8,9,12,13}{blue}
\lvertcols{1,4,5,6,7,8,9,12,13}{blue}
\draw[|-] (0,0)--(0,2);
\draw(0,1)node[left]{$\si_{ij}$};
\end{scope}
\begin{scope}[shift={(0,-21)}]	
\stlines{1/1,4/4,5/5}
\uarcs{6/7,8/9,12/13}
\darcs{6/7,8/9,12/13}
\udotteds{1/4,9/12}
\ldotteds{1/4,9/12}
\uvertcols{1,4,5,6,7,8,9,12,13}{red}
\lvertcols{1,4,5,6,7,8,9,12,13}{blue}
\draw[|-|] (0,0)--(0,2);
\draw(0,1)node[left]{$\al$};
\end{scope}
\begin{scope}[shift={(0,-23)}]	
\stlines{1/1,4/4}
\darcxhalf6{9.5}1
\uarcxhalfr{9.5}{13}1
\uarcs{5/6,7/8}
\uarchalf9{9.5}
\uarchalfr{11.5}{12}
\darcx57{.5}
\darcs{8/9,12/13}
\udotteds{1/4,9/12}
\ldotteds{1/4,9/12}
\uvertcols{1,4,5,6,7,8,9,12,13}{blue}
\lvertcols{1,4,5,6,7,8,9,12,13}{blue}
\draw[|-] (0,0)--(0,2);
\draw(0,1)node[left]{$\be$};
\end{scope}
\begin{scope}[shift={(17,0)}]	
\stlines{1/1,4/4,5/5}
\uarcs{6/7,8/9,10/11,14/15}
\darcs{6/7,8/9,10/11,14/15}
\udotteds{1/4,11/14}
\ldotteds{1/4,11/14}
\uvertcols{1,4,5,6,7,8,9,10,11,14,15}{red}
\lvertcols{1,4,5,6,7,8,9,10,11,14,15}{blue}
\draw[|-|] (0,0)--(0,2);
\draw(0,1)node[left]{$\al$};
\end{scope}
\begin{scope}[shift={(17,-2)}]	
\stlines{1/1,4/4,5/5,6/6,7/8,8/7,9/9,10/10,11/11,14/14,15/15}
\udotteds{1/4,11/14}
\ldotteds{1/4,11/14}
\uvertcols{1,4,5,6,7,8,9,10,11,14,15}{blue}
\lvertcols{1,4,5,6,7,8,9,10,11,14,15}{blue}
\draw[|-] (0,0)--(0,2);
\draw(0,1)node[left]{$\si_{ij}$};
\end{scope}
\begin{scope}[shift={(17,-6)}]	
\stlines{1/1,4/4,5/5}
\uarcs{6/7,8/9,10/11,14/15}
\darcs{6/7,8/9,10/11,14/15}
\udotteds{1/4,11/14}
\ldotteds{1/4,11/14}
\uvertcols{1,4,5,6,7,8,9,10,11,14,15}{red}
\lvertcols{1,4,5,6,7,8,9,10,11,14,15}{blue}
\draw[|-|] (0,0)--(0,2);
\draw(0,1)node[left]{$\al$};
\end{scope}
\begin{scope}[shift={(17,-8)}]	
\stlines{1/1,4/4}
\darcxhalf5{10}1
\uarcxhalfr{10}{15}1
\uarcs{5/6,7/8,9/10}
\uarchalf{11}{11.5}
\uarchalfr{13.5}{14}
\darcs{10/11,14/15}
\darcx68{.3}
\darcx79{.6}
\udotteds{1/4,11/14}
\ldotteds{1/4,11/14}
\uvertcols{1,4,5,6,7,8,9,10,11,14,15}{blue}
\lvertcols{1,4,5,6,7,8,9,10,11,14,15}{blue}
\draw[|-] (0,0)--(0,2);
\draw(0,1)node[left]{$\be$};
\end{scope}
\begin{scope}[shift={(17,-15)}]	
\stlines{1/1,4/4,5/5}
\uarcs{6/7,8/9,12/13}
\darcs{6/7,8/9,12/13}
\udotteds{1/4,9/12}
\ldotteds{1/4,9/12}
\uvertcols{1,4,5,6,7,8,9,12,13}{red}
\lvertcols{1,4,5,6,7,8,9,12,13}{blue}
\draw[|-|] (0,0)--(0,2);
\draw(0,1)node[left]{$\al$};
\end{scope}
\begin{scope}[shift={(17,-17)}]	
\stlines{1/1,4/4,5/5,6/7,7/6,8/8,9/9,12/12,13/13}
\udotteds{1/4,9/12}
\ldotteds{1/4,9/12}
\uvertcols{1,4,5,6,7,8,9,12,13}{blue}
\lvertcols{1,4,5,6,7,8,9,12,13}{blue}
\draw[|-] (0,0)--(0,2);
\draw(0,1)node[left]{$\si_{ij}$};
\end{scope}
\begin{scope}[shift={(17,-21)}]	
\stlines{1/1,4/4,5/5}
\uarcs{6/7,8/9,12/13}
\darcs{6/7,8/9,12/13}
\udotteds{1/4,9/12}
\ldotteds{1/4,9/12}
\uvertcols{1,4,5,6,7,8,9,12,13}{red}
\lvertcols{1,4,5,6,7,8,9,12,13}{blue}
\draw[|-|] (0,0)--(0,2);
\draw(0,1)node[left]{$\al$};
\end{scope}
%
%
%
%
\rectangle{-2}{-24}{33}{4}
\rectangle{-2}{-24}{33}{-11}
\rectangle{-2}{-24}{15}{4}
\draw(6.5,3)node{(i)};
\draw(6.5,-12)node{(ii)};
\draw(24,3)node{(iii)};
\draw(24,-12)node{(iv)};
\end{tikzpicture}
\end{center}
\vspace{-5mm}
\caption{Diagrammatic verification that $\al\si_{ij}=\al\star\be$, where $\be\in\dla \Et(D_r)\dra$, as in the proof of Proposition~\ref{prop:Ir0=<E(Dr)>}; see the text for more details.  In all cases, red vertices are ordered $i_1,\ldots,i_r,a_1,b_1,\ldots,a_s,b_s$, and blue vertices are ordered $j_1,\ldots,j_r,c_1,d_1,\ldots,c_s,d_s$.}
\label{fig:<E(D_r)>}
\end{figure}
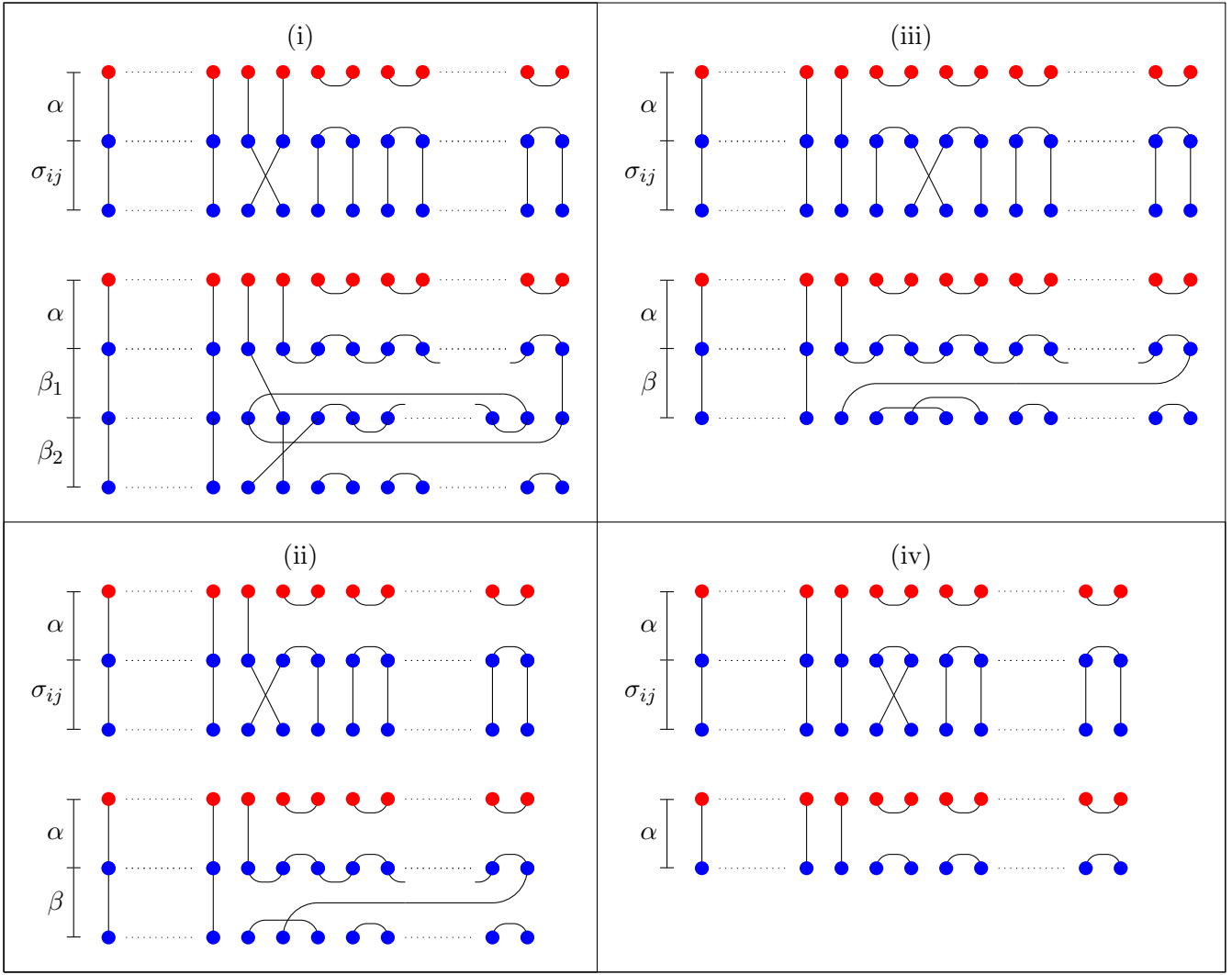

\pf By Proposition \ref{prop:Ir0=<Dr>}, it suffices to show that $D_r\sub\dla\Et(D_r)\dra$.  By Proposition \ref{prop:regular}, $D_r=\DClass r0$ is a regular $\Dt$-class of $\Bnt$, so we may choose an idempotent $\ve\in\Et(D_r)$.  
Since $D_r=D_\ve=\S_n\ve\S_n$, by Proposition \ref{prop:Green_Bn}, it suffices to show that $\lam\ve\rho\in\dla\Et(D_r)\dra$ for all $\lam,\rho\in\S_n$.  In fact, by a simple induction on the length of $\lam$ and $\rho$ as products of transpositions, it suffices to show that
\bit
\item[(I)] for all $\al\in D_r$ and all $\oijn$, $\al\si_{ij}=\al\star\be$ for some $\be\in\dla\Et(D_r)\dra$, and
\item[(II)] for all $\al\in D_r$ and all $\oijn$, $\si_{ij}\al=\be\star\al$ for some $\be\in\dla\Et(D_r)\dra$,
\eit
where we denote by $\si_{ij}\in\S_n$ the transposition that interchanges $i$ and $j$.  By symmetry, it suffices just to prove~(I).  So let $\al\in D_r$ and $\oijn$ be arbitrary, and write $\al$ as in Equation \eqref{eq:al} on Page \pageref{eq:al}.
Recall that $r,s\geq1$.  We now consider four separate cases:
\bmc2
\item[(i)] $i,j\in\codom(\al)$,
\item[(ii)] $i\in\codom(\al)$ but $j\in\bn\sm\codom(\al)$,
\item[(iii)] $i,j\in\bn\sm\codom(\al)$ but $(i,j)\not\in\coker(\al)$,
\item[(iv)] $(i,j)\in\coker(\al)$.
\emc
We show that in all cases, $\al\si_{ij}=\al\star\be$ for some $\be\in\dla\Et(D_r)\dra$.  First, we consider case (i).  Relabelling the vertices, if necessary, we may assume that $(i,j)=(j_{r-1},j_r)$.  In Figure \ref{fig:<E(D_r)>}(i), we show that $\al\si_{ij}=\al\be_1\be_2$ for some $\be_1,\be_2\in\Et(D_r)$ with $\tau(\al,\be_1,\be_2)=0$, giving $\al\si_{ij}=\al\star(\be_1\star\be_2)$, as required (we leave it to the reader to verify that $\be_1,\be_2\in\Et(\Bnt)$).  Similarly, for cases (ii), (iii), (iv), we may assume that $(i,j)=(j_r,a_1)$, $(i,j)=(b_1,a_2)$, $(i,j)=(a_1,b_1)$, respectively.  In Figure~\ref{fig:<E(D_r)>}, we show that $\al\si_{ij}=\al\star\be$ for some $\be\in\Et(D_r)$ in cases (ii) and (iii), and that $\al\si_{ij}=\al$ in case (iv).  As noted above, this completes the proof. \epf

\ms
\begin{rem}
The trick in the above proof, of considering expressions of the form $\al\si_{ij}$ and $\si_{ij}\al$, bears some resemblance to the proof of \cite[Lemma 1.2]{AM2005}.
\end{rem}

The proof of the next result uses several ideas and results from \cite{Gray2008}; see also \cite{GrayHowieIssuePaper}.


\ms
\begin{prop}\label{prop:Bob}
Suppose $r\in I(n)\sm\{0,n\}$.  Then $\Ideal r0$ is idempotent-generated, and
\[
\rank(\Ideal r0)=\idrank(\Ideal r0)=\rho_{nr},
\]
where the numbers $\rho_{nr}$ are defined in Proposition \ref{prop:combinatorics_Bn}.
\end{prop}

\pf For simplicity, write $D=\DClass r0$ and $I=\Ideal r0$.  So $I=\dla \Et(D)\dra$, by Proposition \ref{prop:Ir0=<E(Dr)>}.  The \emph{principal factor} of $D$, denoted $D^\circ$, is the semigroup on the set $D\cup\{0\}$, with multiplication $\circ$ defined, for $\al,\be\in D$, by
\[
\al\circ0=0\circ\al=0\circ0=0 \AND \al\circ\be = \begin{cases}
\al\star\be &\text{if $\al\star\be\in D$}\\
0 &\text{otherwise.}
\end{cases}
\]
Suppose the $\Rt$- and $\Lt$-classes contained in $D$ are $\set{R_j}{j\in J}$ and $\set{L_k}{k\in K}$, where $J\cap K=\emptyset$.  The \emph{Graham-Houghton graph} of $D^\circ$  is the (bipartite) graph $\Delta=\Delta(D^\circ)$ with vertex set $J\cup K$ and edge set $\bigset{\{j,k\}}{\text{$R_j\cap L_k$ contains an idempotent}}$.  We note that $\Delta$ is \emph{balanced}, in the sense that $|J|=|K|$; this common value is equal to $\rho_{nr}$, by Proposition \ref{prop:combinatorics_Bn} and Corollary \ref{cor:Green_Bnt}.  By \cite[Theorem 40]{DEEFHHL1}, each $\Rt$- and $\Lt$-class in~$D$ contains the same number of idempotents; this number was denoted $b_{nr}$ in \cite{DEEFHHL1}, and a recurrence relation was given for these numbers.  It follows that $\Delta$ is \emph{$b_{nr}$-regular}, in the sense that each vertex of $\Delta$ is adjacent to $b_{nr}$ other vertices.  Since $n\geq3$ (as $I(n)\sm\{0,n\}$ is non-empty), we have $b_{nr}\geq2$.  It was shown in \cite[Lemma 3.1]{Gray2008} that being $k$-regular with $k\geq2$ implies that $\Delta$ satisfies the so-called \emph{Strong Hall Condition}:
\[
\text{for all $\emptyset\subsetneq H\subsetneq J$, $|N(H)|>|H|$, where $N(H)$ is the set of all vertices adjacent to a vertex from $H$.}
\]
We also note that $\Delta$ is \emph{connected}; indeed, this follows from the fact that $D^\circ$ is idempotent-generated, as explained in \cite[p61]{Gray2008}.  Since $\Delta$ is connected and balanced and satisfies the Strong Hall Condition, \cite[Lemma 2.11]{Gray2008} gives $\rank(D^\circ)=\idrank(D^\circ)=|J|=|K|=\rho_{nr}$.  But, since $I=\dla D\dra$, it follows that $\rank(I)=\rank(D^\circ)$ and $\idrank(I)=\idrank(D^\circ)$. \epf

Next we consider the ideals $\Ideal rk$ where $r,k>0$.  First we need a technical lemma.


\ms
\begin{lemma}\label{lem:alphabeta2}
Let $\al\in \BnSn$.  
\bit
\itemit{i} If $\rank(\al)>0$, then $\al=\al\star\be$ for some $\be\in D_\al$.
\itemit{ii} If $\rank(\al)=0$, then $\al=\al\star\be$ for some $\be\in D_2$.
\eit
\end{lemma}

\begin{figure}[h]
\begin{center}
\begin{tikzpicture}[scale=.5]
\begin{scope}[shift={(0,0)}]	
\stlines{1/1,4/4,5/5}
\uarcs{6/7,8/9,12/13}
\darcs{6/7,8/9,12/13}
\udotteds{1/4,9/12}
\ldotteds{1/4,9/12}
\uvertcols{1,4,5,6,7,8,9,12,13}{red}
\lvertcols{1,4,5,6,7,8,9,12,13}{blue}
\draw[|-|] (0,0)--(0,2);
\draw(0,1)node[left]{$\al$};
\end{scope}
\begin{scope}[shift={(0,-4)}]	
\stlines{1/1,4/4,5/5}
\uarcs{6/7,8/9,12/13}
\darcs{6/7,8/9,12/13}
\udotteds{1/4,9/12}
\ldotteds{1/4,9/12}
\uvertcols{1,4,5,6,7,8,9,12,13}{red}
\lvertcols{1,4,5,6,7,8,9,12,13}{blue}
\draw[|-|] (0,0)--(0,2);
\draw(0,1)node[left]{$\al$};
\end{scope}
\begin{scope}[shift={(0,-6)}]	
\stlines{1/1,4/4}
\darcxhalf581
\uarcxhalfr8{13}1
\uarcs{5/6,7/8}
\uarchalf9{9.5}
\uarchalfr{11.5}{12}
\darcs{6/7,8/9,12/13}
\udotteds{1/4,9/12}
\ldotteds{1/4,9/12}
\uvertcols{1,4,5,6,7,8,9,12,13}{blue}
\lvertcols{1,4,5,6,7,8,9,12,13}{blue}
\draw[|-] (0,0)--(0,2);
\draw(0,1)node[left]{$\be$};
\end{scope}
\begin{scope}[shift={(20,0)}]	
\uarcs{1/2,3/4,7/8}
\darcs{1/2,3/4,7/8}
\udotteds{4/7}
\ldotteds{4/7}
\uvertcols{1,2,3,4,7,8}{red}
\lvertcols{1,2,3,4,7,8}{blue}
\draw[|-|] (0,0)--(0,2);
\draw(0,1)node[left]{$\al$};
\end{scope}
\begin{scope}[shift={(20,-4)}]	
\uarcs{1/2,3/4,7/8}
\darcs{1/2,3/4,7/8}
\udotteds{4/7}
\ldotteds{4/7}
\uvertcols{1,2,3,4,7,8}{red}
\lvertcols{1,2,3,4,7,8}{blue}
\draw[|-|] (0,0)--(0,2);
\draw(0,1)node[left]{$\al$};
\end{scope}
\begin{scope}[shift={(20,-6)}]
\stlines{1/1}
\darcxhalf251
\uarcxhalfr581
\uarcs{2/3}
\uarchalf4{4.5}
\uarchalfr{6.5}{7}
\darcs{3/4,7/8}
\udotteds{4/7}
\ldotteds{4/7}
\uvertcols{1,2,3,4,7,8}{blue}
\lvertcols{1,2,3,4,7,8}{blue}
\draw[|-] (0,0)--(0,2);
\draw(0,1)node[left]{$\be$};
\end{scope}
\end{tikzpicture}
\end{center}
\vspace{-5mm}
\caption{Diagrammatic verification that $\al=\al\be$ from the proof of Lemma \ref{lem:alphabeta2} ($r>0$ on the left, $r=0$ on the right); see the text for more details.  Red vertices are ordered $i_1,\ldots,i_r,a_1,b_1,\ldots,a_s,b_s$, and blue vertices are ordered $j_1,\ldots,j_r,c_1,d_1,\ldots,c_s,d_s$.}
\label{fig:alphabeta2}
\end{figure}
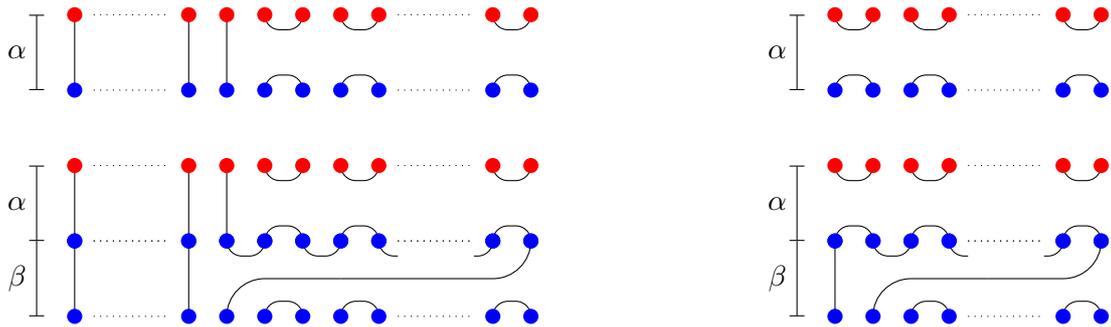

\pf Write $\al$ as in Equation \eqref{eq:al} on Page \pageref{eq:al}.
In Figure \ref{fig:alphabeta2}, we demonstrate the existence of $\be$ (of the desired rank) such that $\al=\al\be$ with $\tau(\al,\be)=0$. \epf


\newpage

\ms
\begin{prop}\label{prop:Mrk}
Let $r\in I(n)\sm\{0\}$, and let $k\geq1$.  Put
\[
\Minimal rk=\bigcup_{s\in I(r)\atop k\leq l< 2k} \DClass sl = \set{(l,\al)\in\Bnt}{k\leq l<2k,\ \rank(\al)\leq r}.
\]
Then
\bit
\itemit{i} $\Ideal rk = \dla\Minimal rk\dra$,
\itemit{ii} any generating set for $\Ideal rk$ contains $\Minimal rk$, so $\Minimal rk$ is the unique minimal (with respect to size or inclusion) generating set for $\Ideal rk$,
\itemit{iii} $\ds{\rank(\Ideal rk)=|\Minimal rk|=k\cdot\sum_{s\in I(r)}\de_{ns}}$, where the numbers $\de_{ns}$ are defined in Proposition \ref{prop:combinatorics_Bn}.
\eit
\end{prop}

\pf We begin with (i).  Let $(i,\al)\in\Ideal rk$ be arbitrary.  If $k\leq i<2k$, then $(i,\al)\in\Minimal rk$, so suppose $i\geq2k$.  Write $i=qk+s$, where $q\in\N$ and $k\leq s<2k$.  By Lemma \ref{lem:alphabeta2}, there exists $\be\in I_r$ such that $\al=\al\be$ with $\tau(\al,\be)=0$.  But then
\[
(i,\al) = (s,\al)\star \underbrace{(k,\be)\star\cdots\star(k,\be)}_q \in\dla\Minimal rk\dra.
\]
This completes the proof of (i).
For (ii), suppose $\Ga$ is an arbitrary generating set for $\Ideal rk$.  Let $(i,\al)\in\Minimal rk$ be arbitrary, and consider an expression
\[
(i,\al)=(i_1,\al_1)\star\cdots\star(i_t,\al_t)=(i_1+\cdots+i_t+\tau(\al_1,\ldots,\al_t),\al_1\cdots\al_t),
\]
where $(i_1,\al_1),\ldots,(i_t,\al_t)\in\Ga$.  Since $i_1,\ldots,i_t\geq k$ and since $i<2k$, it follows that $t=1$, so that $(i,\al)=(i_1,\al_1)\in\Ga$, giving (ii).  It follows immediately from (i) and (ii) that $\rank(\Ideal rk)=|\Minimal rk|$.  The formula for $|\Minimal rk|$ follows from the fact that $|\DClass sl|=|D_s|=\de_{ns}$ (see Corollary \ref{cor:Green_Bnt} and Proposition \ref{prop:combinatorics_Bn}). \epf

\ms
\begin{rem}
The generating set $\Minimal52$ of the ideal $\Ideal52$ of $\B_7^\tau$ is pictured in Figure \ref{fig:partial_order}.
\end{rem}

\ms
\begin{prop}\label{prop:M0k}
Suppose $n$ is even, and let $k\in\N$ be arbitrary.  Put $\Minimal0k=\DClass0k\cup\cdots\cup\DClass0{2k}$.  Then
\bit
\itemit{i} $\Ideal 0k = \dla\Minimal 0k\dra$,
\itemit{ii} any generating set for $\Ideal 0k$ contains $\Minimal 0k$, so $\Minimal 0k$ is the unique minimal (with respect to size or inclusion) generating set for $\Ideal 0k$,
\itemit{iii} $\rank(\Ideal 0k)=|\Minimal 0k|=(k+1)\cdot\de_{n0}$, where the numbers $\de_{n0}$ are defined in Proposition \ref{prop:combinatorics_Bn}.
\eit
\end{prop}

\pf We omit the proof as it is very similar to that of Proposition \ref{prop:Mrk}.  The main difference is that we apply Lemma \ref{lem:alphabeta} instead of Lemma \ref{lem:alphabeta2}.  This explains the factor of $k+1$ in the expression for $\rank(\Ideal0k)$.~\epf

\ms
\begin{rem}
Note that $\Minimal00=\DClass00=D_0$.  We saw that $\Ideal00=\dla D_0\dra$ in Proposition \ref{prop:Ir0=<Dr>}.
\end{rem}

For convenience, we gather the results on ranks of principal ideals into a single theorem.

\ms
\begin{thm}\label{thm:rank}
Let $n\geq3$, $r\in I(n)$ and $k\in\N$.  Then
\[
\rank(\Ideal rk) = \begin{cases}
4 &\text{if $r=n$ and $k=0$}\\
\rho_{nr} &\text{if $0<r<n$ and $k=0$}\\
(k+1)\cdot\de_{n0} &\text{if $r=0$}\\
k\cdot\sum_{s\in I(r)}\de_{ns} &\text{if $r>0$ and $k>0$,}
\end{cases}
\]
where the numbers $\rho_{nr},\de_{nr}$ are defined in Proposition \ref{prop:combinatorics_Bn}. Further, $\Ideal rk$ is idempotent-generated if and only if $0<r<n$ and $k=0$, in which case $\idrank(\Ideal rk)=\rank(\Ideal rk)$.  \epf
\end{thm}

\ms
\begin{rem}
An obvious necessary condition for an ideal $I$ of an arbitrary semigroup $S$ to be idempotent-generated is that there must be idempotents in any maximal $\J$-class of $I$.  Since idempotents of $\Bnt$ can only exist in $\Jt=\Dt$-classes of the form $D_{r;0}$, it follows that any idempotent-generated ideal of $\Bnt$ is a principal ideal (of the form described in Theorem \ref{thm:rank}).
\end{rem}

\subsection{Applications}\label{sect:applications}

A famous result of Howie \cite{Howie1966} states that the singular ideal $\TnSn$ of the full transformation semigroup $\T_n$ is idempotent-generated.  In fact, the idempotent-generated subsemigroup $\la E(\T_n)\ra$ is equal to $\{1\}\cup(\TnSn)$.  This is true also of the Brauer monoid $\B_n$: specifically, $\la E(\B_n)\ra=\{1\}\cup(\BnSn)$, as shown in \cite{Maltcev2007}, where a presentation for $\BnSn$ was also given.  Similar results for other diagram semigroups appear in \cite{JEpnsn,EastGray,DEG}.

We now apply the results of previous sections to explore the analogous situation for the twisted Brauer monoid $\Bnt$.  This is more complicated, and it is not the case that the singular ideal $\BntSn$ is idempotent-generated.  We may still calculate the rank of this singular ideal, and we also describe the idempotent-generated subsemigroup $\dla\Et(\Bnt)\dra$, and calculate its rank and idempotent rank (which are equal).  We also deduce the above-mentioned result that $\BnSn$ is idempotent-generated.

\ms
\begin{thm}\label{thm:SingBnt}
If $n\geq3$, then $\rank(\BntSn)={n\choose2}+n!$.
\end{thm}

\pf Note that $\BntSn=\Ideal{n-2}0\cup\Ideal n1$.  By (the proof of) Proposition \ref{prop:Bob}, we may choose 
a subset $\Si\sub\DClass{n-2}0$ with $\Ideal{n-2}0=\dla\Si\dra$ and $|\Si|=\rank(\Ideal{n-2}0)=\rho_{n,n-2}={n\choose2}$.  Now put $\Ga=\Si\cup\DClass n1$.  Since $\dla\Si\dra=\Ideal{n-2}0\supseteq \DClass z1\cup\cdots\cup\DClass{n-2}1$, it follows that $\dla\Ga\dra\supseteq\Minimal n1$, so $\dla\Ga\dra\supseteq\Ideal n1$.  Thus, $\BntSn=\Ideal{n-2}0\cup\Ideal n1=\dla\Ga\dra$.  In particular,
\begin{equation}\label{eq0}
\rank(\BntSn)\leq|\Ga|=|\Si|+|\DClass n1|={n\choose2}+n!.
\end{equation}
Conversely, suppose $\Xi$ is an arbitrary generating set for $\BntSn$.  Let $\al\in\DClass{n-2}0$ be arbitrary, and consider an expression
\[
\al=(0,\al) = (i_1,\al_1)\star\cdots\star(i_k,\al_k)=(i_1+\cdots+i_k+\tau(\al_1,\ldots,\al_k),\al_1\cdots\al_k),
\]
where $(i_1,\al_1),\ldots,(i_k,\al_k)\in\Xi$.  Then we must have
\[
i_1=\cdots=i_k=\tau(\al_1,\ldots,\al_k)=0 \AND \al_1\cdots\al_k=\al.
\]
Then for any $j\in\bk$, $n-2=\rank(\al)=\rank(\al_1\cdots\al_k)\leq\rank(\al_j)\leq n-2$.  In particular, $(i_j,\al_j)\in\DClass{n-2}0$ for each $j\in\bk$.  We have shown that $\DClass{n-2}0\sub\dla\Xi\cap\DClass{n-2}0\dra$.  It follows that $\Ideal{n-2}0=\dla\DClass{n-2}0\dra\sub\dla\Xi\cap\DClass{n-2}0\dra$.
In particular,
\begin{equation}\label{eq1}
|\Xi\cap\DClass{n-2}0| \geq \rank(\Ideal{n-2}0)={n\choose2}.
\end{equation}
Next, let $\si\in\S_n$ be arbitrary.  As in the proof of Proposition \ref{prop:Mrk}, consideration of an expression for $(1,\si)$ as a product of elements from $\Xi$ shows that, in fact, $(1,\si)\in\Xi$.  In particular, it follows that $\DClass n1\sub\Xi$, so
\begin{equation}\label{eq2}
|\Xi\sm\DClass{n-2}0| \geq|\DClass n1|=n!.
\end{equation}
Adding \eqref{eq1} and \eqref{eq2}, we obtain $|\Xi|\geq{n\choose2}+n!$.  Since $\Xi$ was an arbitrary generating set for $\BntSn$, it follows that $\rank(\BntSn)\geq{n\choose2}+n!$.  Combined with \eqref{eq0}, this completes the proof. \epf


We now describe the idempotent-generated subsemigroup of $\Bnt$ and also derive a formula for its rank and idempotent rank.

\ms
\begin{thm}\label{thm:IGBnt}
Let $n\geq3$ and let $S=\dla\Et(\Bnt)\dra$ be the idempotent-generated subsemigroup of $\Bnt$.  Then
\[
S = \{1\}\cup\Ideal{n-2}0 = \{1\} \cup (\N\times(\BnSn)) = \{1\}\cup\set{(i,\al)}{i\in\N,\ \al\in\BnSn},
\]
and $\rank(S)=\idrank(S)={n\choose2}+1$.
\end{thm}

\pf Since $1\in\Et(\Bnt)$, and since $\Ideal{n-2}0$ is idempotent-generated, by Proposition \ref{prop:Ir0=<E(Dr)>}, it is clear that $\{1\}\cup\Ideal{n-2}0\sub S$.  To show the reverse containment, it suffices to show that $S\sm\Ideal{n-2}0=\{1\}$.  So suppose $(i,\al)\in S\sm\Ideal{n-2}0$.  In particular, $\al\in D_n=\S_n$, and we have
\[
(i,\al) = \al_1\star\cdots\star\al_k = (\tau(\al_1,\ldots,\al_k),\al_1\cdots\al_k)
\]
for some idempotents $\al_1,\ldots,\al_k\in\Et(\Bnt)$.  
(Recall that $\Et(\Bnt)\sub E(\B_n)$.)  Since $\al_1\cdots\al_k=\al\in\S_n$, and since $\BnSn$ is an ideal of $\B_n$, it follows that $\al_1,\ldots,\al_k\in\S_n$.  In particular, $\tau(\al_1,\ldots,\al_k)=0$.  But also $E(\S_n)=\{1\}$, as $\S_n$ is a group.  So $\al_1=\cdots=\al_k=1$, and $(i,\al)= (\tau(\al_1,\ldots,\al_k),\al_1\cdots\al_k)=(0,1)=1$, as required.
The statement about the rank and idempotent rank follows immediately from Proposition \ref{prop:Bob} and the obvious fact that $\Ideal{n-2}0=\dla\Si\dra \iff S=\{1\}\cup\Ideal{n-2}0=\dla\{1\}\cup \Si\dra$. \epf

As a final application, we prove the following result, which is a (slight) strengthening of a result from \cite{Maltcev2007}.

\ms
\begin{thm}[cf.~Maltcev and Mazorchuk \cite{Maltcev2007}]\label{thm:MM}
Let $n\geq3$.  The singular part $\BnSn$ of the Brauer monoid~$\B_n$ is idempotent-generated.  In fact, 
\[
\{1\}\cup(\BnSn) = \la E(\B_n)\ra = \la \Et(\B_n)\ra.
\]
\end{thm}

\pf 
First, it is clear that $\la \Et(\B_n)\ra\sub \la E(\B_n)\ra$.
Next note that, since no non-identity element of $\S_n$ is a product of idempotents (from $\B_n$), we have $\la E(\B_n)\ra\sub\{1\}\cup(\BnSn)$.  Finally, suppose $\al\in\BnSn$ is arbitrary.  Then, by Theorem \ref{thm:IGBnt},
$
(0,\al)=\al = \al_1\star\cdots\star\al_k=(\tau(\al_1,\ldots,\al_k),\al_1\cdots\al_k)
$
for some idempotents $\al_1,\ldots,\al_k\in\Et(\Bnt)$.  In particular, $\al=\al_1\cdots\al_k\in\la \Et(\B_n)\ra$.  (And also $\tau(\al_1,\ldots,\al_k)=0$.)  \epf

\ms
\begin{rem}\label{rem:MM1}
One may easily check that $\la E(\B_2)\ra=\{1\}\cup(\B_2\sm\S_2)\not=\{1\}=\la\Et(\B_2)\ra$.  In the above proof, we showed that every element of $\BnSn$ (with $n\geq3$) may be written as a product of idempotents from $\Et(\B_n)$ in such a way that no floating components are created in the formation of the product.  We note that this also follows from \cite[Proposition 2]{Maltcev2007} or \cite[Proposition 8.7]{EastGray}, but we omit the details.
%
\end{rem}




\subsection*{Acknowledgement}

We kindly thank the anonymous referee for his/her careful reading of the article, and for several suggestions leading to greater readability.

\footnotesize
\def\bibspacing{-1.1pt}
\bibliography{biblio}
\bibliographystyle{plain}
\end{document}